\documentclass[12pt]{amsart}
\usepackage{amssymb}
\usepackage{amsmath}
\usepackage{pictex}
\newcommand\1{\mathbf{1}}

\newcommand\al{\alpha}
\newcommand\AND{\quad\mbox{and}\quad}
\newcommand\BB{\mathcal B}
\newcommand\bd{\partial}
\newcommand\be{\beta}

\newcommand\CC{\mathcal C}
\newcommand\cf{\curlywedge}

\newcommand\de{\delta}
\newcommand\DL{\mathsf{DL}}
\newcommand\dn{\mathfrak{d}}
\newcommand\dps{\displaystyle}
\newcommand\Ex{\mathsf{E}}

\newcommand\geo[1]{\overline{#1}}

\newcommand\HH{{\mathcal H}}
\newcommand\hor{\mathfrak{h}}

\newcommand\la{\lambda}
\newcommand\lle{\preccurlyeq}
\newcommand\lra{\leftrightarrow}

   \newcommand\MM{\mathcal M}
\newcommand\om{\omega}
\newcommand\mm{\mathsf{m}}
\newcommand\nn{\mathfrak{n}}
\newcommand\Prb{\operatorname{\sf Pr}}
   \newcommand\R{\mathbb R}
\newcommand\Rest{\operatorname{\sf Rest}}
\newcommand\scs{\scriptstyle}

\newcommand\spn{\mathfrak{s}}
\newcommand\supp{\operatorname{\sf supp}}
\newcommand\T{\mathbb T}

\newcommand\tm{\mathbf{t}}
\newcommand\up{\mathfrak{u}}
\newcommand\wh{\widehat}
\newcommand\wt{\widetilde}
\newcommand\Z{\mathbb Z}

\numberwithin{equation}{section}

\newtheoremstyle{mythm}
  {9pt}
  {9pt}
  {\itshape}
  {0pt}
  {\bfseries}
  {}
  { }
  {\thmnumber{(#2)}\thmname{ #1}\thmnote{ #3}}

\newtheoremstyle{mydef}
  {9pt}
  {9pt}
  {\normalfont}
  {0pt}
  {\bfseries}
  {}
  { }
  {\thmnumber{(#2)}\thmname{ #1}\thmnote{ #3}}

\theoremstyle{mythm}
\newtheorem{thm}[equation]{Theorem.}
\newtheorem{pro}[equation]{Proposition.}
\newtheorem{lem}[equation]{Lemma.}
\newtheorem{cor}[equation]{Corollary.}

\theoremstyle{mydef}
\newtheorem{dfn}[equation]{Definition.}

\begin{document}
\title{\large Green kernel estimates and the full Martin boundary 
for random walks on lamplighter groups and Diestel-Leader graphs}
\author{\bf Sara BROFFERIO and Wolfgang WOESS}
\address{\parbox{1.4\linewidth}{Institut f\"ur Mathematik C, 
Technische Universit\"at Graz\\
Steyrergasse 30, A-8010 Graz, Austria\\ \ \\
\it E-mail: \tt brofferio@weyl.math.tu-graz.ac.at\\ 
\hspace*{1.12cm} woess@weyl.math.tu-graz.ac.at
}}
\date{29 February 2004}
\thanks{Supported by European Commission, Marie Curie Fellowship
HPMF-CT-2002-02137 and partially by FWF (Austrian Science Fund)
project P15577}
\begin{abstract}
We determine the precise asymptotic behaviour (in space) of the Green kernel of 
simple random walk with drift on the Diestel-Leader graph $\DL(q,r)$,
where $q,r \ge 2$. 
The latter is the horocyclic product of two homogeneous trees with 
respective degrees $q+1$ and $r+1$. 
When $q=r$, it is the Cayley graph of the wreath product 
(lamplighter group) $\Z_q \wr \Z$ with respect to a natural set of 
generators. 
We describe the full Martin compactification of these random walks
on $\DL$-graphs and, in particular, lamplighter groups. 
This completes and provides a better approach to previous results
of {\sc Woess,} who has determined all minimal positive harmonic functions.
\end{abstract}
\maketitle

\markboth{\sf S. Brofferio and W. Woess}{\sf Green kernel estimates}
\baselineskip 15pt

\section{Introduction}\label{intro}

Consider the additive group $\Z$ of all integers as a two-way-infinite
road where at each point there is a lamp that may be switched on in one of 
$q$ different intensities (states) $\in \Z_q  = \{ 0, \dots, q-1 \}$, the 
group of integers modulo $q$. At the beginning, all lamps are in state $0$ 
(switched off), and a lamplighter starts at some point of $\Z$. He chooses at 
random among the following actions (or a suitable combination thereof): 
he can move to a neighbour point in $\Z$, or he can change the intensity 
of the lamp at the actual site to a different state.
When the process evolves, we have to keep track of the position 
$k \in \Z$ of the lamplighter and of the finitely supported configuration
$\eta: \Z \to \Z_q$ that describes the states of all lamps.
The set $\Z_q \wr \Z$ of all pairs $(\eta, k)$ of this type
carries the structure of a semidirect product of $\Z$ with the
additive group $\CC$ of all configurations, on which $\Z$ acts by
translations. This is often called the \emph{lamplighter group};
the underlying algebraic construction is the \emph{wreath product} of two
groups. 


Random walks on lamplighter groups have been a well-studied subject
in recent years, see {\sc Kaimanovich and Vershik} \cite{KaiVer} and 
{\sc Kaimanovich} \cite{Kai} (Poisson boundary $\equiv$ bounded harmonic
functions), {\sc Lyons, Pemantle and Peres} 
\cite{LyoPemPer}, {\sc Erschler} \cite{Ers}, \mbox{\sc Revelle} \cite{Rev1},
{\sc Bertacchi} \cite{Ber} (rate of escape), {\sc Grigorchuk and \.Zuk},
\cite{GriZuk}, {\sc Dicks and Schick} \cite{DiSc}, {\sc Bartholdi and Woess} 
\cite{BaWo} (spectral theory), {\sc Saloff-Coste and
Pittet} \cite{PitSal1}, \cite{PitSal2}, {\sc Revelle} \cite{Rev2}
(asymptotic behaviour of transition probabilites), and {\sc Woess} \cite{Woe}
(positive harmonic functions).

Here, we shall deal with Green kernel asymptotics and 
positive harmonic functions. Let us briefly outline in general how this
is linked with \emph{Martin boundary theory} of Markov chains.
Consider an arbitrary infinite (connected, locally finite) graph $X$ (e.g. a
Cayley graph of a finitely generated group) and the stochastic 
transition matrix $P = \bigl( p(x,y) \bigr)_{x,y \in X}$  of a random walk 
$Z_n$ on $X$. That is, $Z_n$ is an 
$X$-valued random variable, the position of the random walker at
time $n$, subject to the Markovian transition rule
$$
\Prb[Z_{n+1} = y \mid Z_n=x ] = p(x,y)\,.
$$
The $n$-step transition probabilitiy
$$
p^{(n)}(x,y) = \Prb[Z_n=y \mid Z_0 = x]\,,\quad x, y \in X\,,
$$
is the $(x,y)$-entry of the matrix power $P^n$, with $P^0=I$, the
identity matrix. The \emph{Green kernel} is
$$
G(x,y) = \sum_{n=0}^{\infty} p^{(n)}(x,y)\,, \quad x,y \in X\,.
$$
This is the expected number of visits in the point $y$, when the
random walk starts at $x$. We always consider random walks that
are \emph{irreducible} and \emph{transient}, which amounts to
$$
0 < G(x,y) < \infty \quad \mbox{for all}\; x,y \in X\,.
$$
\emph{Renewal theory} in a wide sense consists in describing the
asymptotic behaviour in space of $G(x,y)$, when $x$ is fixed and
$y$ tends to infinity (or dually, $y$ is fixed and $x$ tends to infinity). 
If we fix a reference point $o \in X$, then the \emph{Martin kernel} is
$$
K(x,y) = G(x,y)/G(o,y)\,,\quad x,y \in X\,.
$$
If we have precise asymptotic estimates in space of the Green kernel, then we 
can also determine the \emph{Martin compactification.} 
This is the smallest metrizable compactification of $X$ containing $X$ 
as a discrete, dense subset, and to which
all functions $K(x,\cdot)$, $x \in X$, extend continuously. The 
\emph{Martin boundary} $\MM = \MM(P)$ is the ideal boundary added to $X$
in this compactification. Thus, $\MM$ consists of the ``directions 
of convergence'' of $K(x,y)$, when $y \to \infty$.
Its significance is that it leads to a complete understanding of the
cone $\HH^+ = \HH^+(P)$ of \emph{positive harmonic functions.}
A function $h: X \to \R$ is called \emph{harmonic,} or \emph{$P$-harmonic,}
if
$$
h = Ph\,,\quad\mbox{where}\quad Ph(x) = \sum_y p(x,y)\,h(y)\,.
$$
A function $h \in \HH^+$ is called \emph{minimal} 
if
$$
h(o) = 1 \AND h \ge h_1 \in \HH^+ \Longrightarrow h_1/h \equiv \text{constant.}
$$
The minimal harmonic functions are the extreme points of the convex base
$\BB = \{ h \in \HH^+ : h(o)=1 \}$ of the cone $\HH^+$.  

The reader is referred to the excellent introduction to Martin boundary theory
by {\sc Dynkin} \cite{Dyn}. A main result of this theory is that
\begin{itemize}
\item Every minimal harmonic function is of the form $K(\cdot,\xi)$, where
$\xi \in \MM$, and the set 
$\MM_{\min} = \{\xi \in \MM : K(\cdot,\xi) \;\mbox{is minimal}\,\}$
is a Borel subset of $\MM\,$;
\item For every $h \in \HH^+$ there is a unique Borel measure $\nu^h$
on $\MM$ such that 
$$
\nu^h(\MM \setminus \MM_{\min}) = 0 \AND 
h(\cdot) = \int_{\MM} K(\cdot,\xi)\,d\nu^h(\xi)\,.
$$
\end{itemize}
The above is an abstract construction of the Martin compactification.
The kind of approach that we have in mind here is the following. The transition
matrix $P$ is adapted to the graph structure, and we want to understand
and describe the Martin compactification in terms of the specific 
geometry of $X$. Results of this type for random walks on various 
classes of graphs and groups, along with many references,
are presented in the book by {\sc Woess} \cite{Wbook}. 

Returning to lamplighter walks, this spirit requires as the starting point
a good understanding of the \emph{geometry} of the wreath product
$\Z_q \wr \Z$ in terms of a suitable Cayley graph of that group.
This is the Diestel-Leader graph $\DL(q,q)$, a special case
of the Diestel-Leader graphs $\DL(q,r)$ ($q,r \ge 2)$, which are explained 
in detail \S \ref{geometry}. Briefly speaking, $\DL(q,r)$
is the \emph{horocyclic product} of the homogeneous trees $\T_q$ and $\T_r$
with degrees $q+1$ and $r+1$, respectively. It is precisely this geometric 
realization of the lamplighter groups in terms of relatively simple objects
such as trees, that allows us to perform many computations.

The random walk with
transition matrix $P_{\al}$ on $\DL(q,r)$ that we consider is the
\emph{simple random walk} (SRW) with an additional drift parameter 
$\al \in (0\,,\,1)$. 
If $r=q$ (the case of the lamplighter group), then this random walk can be 
interpreted  in lamplighter terms as follows. Think of the lamps not placed 
at each  vertex of the  two-way-infinite path $\Z$, but at the middle of each 
edge. Suppose the actual position of the lamplighter is $k \in \Z$.
He first tosses a coin. If ``head'' comes up, which happens with
probability $\al$, he moves to $k+1$ and switches the lamp on the 
transversed edge to a state chosen at random in $\Z_q$.
Otherwise, he moves to $k-1$ and also switches the lamp on the 
transversed edge to a random state.
 
Even when $q \ne r$, the random walk $P_{\al}$ on $\DL(q,r)$
may be interpreted as a lamplighter walk in an extended sense.
Imagine that on each edge of $\Z$, there is a green lamp with $q$ possible
intensities (including ``off'') in $\Z_q$ \emph{plus} a red lamp with 
$r$ possible intensities (including ``off'') in $\Z_r$. 
The rule is that only finitely many lamps 
may be switched on, and in addition, if the lamplighter stands at
$k$, then all red lamps between $k$ and $-\infty$ have to be switched off, 
while all green lamps between $k$ and $+\infty$ must be be switched off.
The lamplighter tosses his $\al$-coin. If ``head'' comes up,
he moves from $k$ to $k+1$ and switches the \emph{green} lamp on 
the transversed edge to a random  state, while switching off the 
red lamp on that edge. Otherwise, he moves to $k-1$ 
and switches the \emph{red} lamp on the transversed edge to a random state, 
while switching off the green lamp sitting there.

Then the random walk $P=P_{\al}$ (whose definition is formalized in 
\eqref{random-walk}) is irreducible and transient. Via our geometric interpretation,
we see that it has natural projections
$P_1 = P_{\al,q}$ and $P_2 =P_{1-\al,r}$ on the two trees used to make up the
graph, and also $\wt P = \wt P_{\al}$ on $\Z$, which describes jsut
the moves of the lamplighter. A good understanding of these projected
walks is crucial for our
approach, and in \S \ref{trees}, we quickly review 
the necessary facts concerning those random walks on $\T_q$ (and $\T_r$).

In \S \ref{renewal}, we derive our main results concerning the asymptotic
behaviour of the Green kernel associated with $P_{\al}$, subsumed in
Theorem \ref{estimates}. The assymptotics are different along different
directions of going to infinity. Also, the drift-free case ($\al = 1/2$)
is substantially different from the other cases ($\al \ne 1/2$).

These results are used in \S \ref{martin} to describe the full
Martin compactification. In the drift-free case, this is the ``natural''
geometric compactification in terms of the two underlying trees. 
Namely, $\DL(q,r)$ is a subgraph of $\T_q \times \T_r$, and the
Martin compactification is the closure of $\DL(q,r)$ in 
$\wh\T_q \times \wh\T_r$, where $\wh\T_q$ and $\wh\T_r$ are the well-known
end-compactifications of the respective trees. However, when $\al \ne 1/2$,
the Martin compactification is larger than ($\equiv$ surjects non-trivially
onto) the geometric one.
The minimal Martin boundary, previously described in \cite{Woe} without
elaborating the directions of convergence, is recovered.

These results can also be adapted to obtain the Martin compactification
for all positive \emph{$t$-harmonic} functions, that satisfy $Ph = t\cdot h$,
where $t \ge \rho(P_{\al}) = \limsup_n p^{(n)}(x,x)^{1/n}$, 
the ``bottom of the positive spectrum''.
The picture at the bottom is that of the drift-free case, while the case 
$t > \rho(P_{\al})$ corresponds to non-vanishing drift. See \S \ref{t-harmonic}.

In the short \S \ref{Harnack}, we present another little by-product of our 
Green kernel estimates, namely we illustrate their
use for showing explicitly that SRW on $\DL(q,q)$ (the lamplighter group)
does not satisfy the \emph{elliptic Harnack inequality}.


In conclusion, let us remark that in general it is significantly 
harder to determine the
whole Martin compactification than to determine the minimal harmonic functions
associated with a random walk, since the former contains more detailed 
analytic-geometric informations than the latter, whose computation 
often has rather an algebraic than an analytic flavour.
Let us also remark that our results provide the first case where 
one can successfully determine the whole Martin
compactifiction of a class of random walks on finitely generated groups
that are solvable, but do not have polynomial growth.

\section{The geometry of Diestel-Leader graphs and lamplighter groups}\label{geometry}

We now explain the structure of the DL-graphs and their relation with the 
wreath products $\Z_q \wr \Z$. 
This section is a short and slightly modified version of \S 2 in \cite{Woe},
included here for the sake of completeness. 

Let $\T = \T_q$ be the homogeneous tree with degree $q+1$, $q \ge 2$.
A \emph{geodesic path}, respectively \emph{geodesic ray}, respectively 
\emph{infinite geodesic} in $\T$ is a finite, respectively one-sided infinite, 
respectively doubly infinite sequence $(x_n)$ of vertices
of $\T$ such that $d(x_i,x_j) = |i-j|$ for all $i, j$, 
where $d(\cdot,\cdot)$ denotes the graph distance. 

Two rays are \emph{equivalent} if their symmetric difference is finite.
An \emph{end} of $\T$ is an equivalence class of rays. The space of 
ends is denoted $\bd \T$, and we write $\wh \T = \T \cup \bd \T$. 
For all $w, z \in \wh \T$, $w \ne z$, there is a unique geodesic $\geo{w\,z}$ 
that connects the two. In particular, if $x \in \T$ and $\xi \in \bd \T$ then 
$\geo{x\,\xi}$ is the ray that starts at $x$ and represents $\xi$.

For $x,y \in \T$, $x \ne y$, we define the \emph{cone}
$\wh \T(x,y) = \{ w \in \wh \T : y \in \geo{x\,w} \}$.
The collection of all cones is the basis of  a topology which 
makes $\wh \T$ a compact, totally disconnected Hausdorff space
with $\T$ as a dense, discrete subset.

We fix a root  $o \in \T$.
If $w, z \in \wh \T$, then their \emph{confluent} $c=w \wedge z$ with
respect to the root vertex $o$ is defined by
$\geo{o\,w} \cap \geo{o\,z} = \geo{o\,c}$. 
Similarly, we choose and fix a \emph{reference end} $\om \in \bd \T$. For 
$z, v \in \wh \T \setminus \{ \om \}$, their confluent $b = v \cf z$ 
with respect to $\om$ is defined by
$\geo{v\,\om} \cap \geo{z\,\om} = \geo{b\,\om}$. We write
$$
z \lle v \quad\mbox{if}\quad z \cf v = z\,.
$$ 
For $x,y \in \T$, we describe their
relative position by the two numbers
$$
\up(x,y) = d(x,x \cf y) \AND \dn(x,y) = d(y,x \cf y)\,.
$$
Thus, $\dn(x,y) = \up(y,x)$. 
In Figure 1, $\up(x,y)$ and $\dn(x,y)$ correspond to the numbers of steps one 
has to take upwards (in direction of $\om$), respectively downwards, 
on the geodesic path from $x$ to $y$. We have $d(x,y) = \up(x,y)+\dn(x,y)$. 

The \emph{Busemann function} $\hor: \T \to \Z$ and the \emph{horocycles} $H_k$
with respect to $\om$ are
$$
\hor(x) = \dn(o,x) - \up(o,x) \AND H_k = \{ x \in \T : \hor(x) = k \}\,.
$$
Every horocycle is infinite. We write $H(x)=H_k$ if $x \in H_k$. 
Every vertex $x$ in $H_k$ has one neighbour 
$x^-$ (its predecessor) in $H_{k-1}$ and $q$ neighbours (its successors)
in $H_{k+1}$. Thus $\lle$ is the transitive closure of the predecessor relation.
We set $\bd^* \T = \bd \T \setminus \{\om\}$. 
$$
\beginpicture 

\setcoordinatesystem units <.7mm,1.04mm>

\setplotarea x from -10 to 104, y from 4 to 84

\arrow <6pt> [.2,.67] from 2 2 to 80 80

\plot 32 32 62 2 /

 \plot 16 16 30 2 /

 \plot 48 16 34 2 /

 \plot 8 8 14 2 /

 \plot 24 8 18 2 /

 \plot 40 8 46 2 /

 \plot 56 8 50 2 /

 \plot 4 4 6 2 /

 \plot 12 4 10 2 /

 \plot 20 4 22 2 /

 \plot 28 4 26 2 /

 \plot 36 4 38 2 /

 \plot 44 4 42 2 /

 \plot 52 4 54 2 /

 \plot 60 4 58 2 /



 \plot 99 29 64 64 /

 \plot 66 2 96 32 /

 \plot 70 2 68 4 /

 \plot 74 2 76 4 /

 \plot 78 2 72 8 /

 \plot 82 2 88 8 /

 \plot 86 2 84 4 /

 \plot 90 2 92 4 /

 \plot 94 2 80 16 /


\setdots <3pt>
\putrule from -4.8 4 to 102 4
\putrule from -4.5 8 to 102 8
\putrule from -2 16 to 102 16
\putrule from -1.7 32 to 102 32
\putrule from -1.7 64 to 102 64
\setdashes <2pt>
\putrule from -1.7 -8 to 102 -8

\put {$\vdots$} at 32 -2
\put {$\vdots$} at 64 -2

\put {$\dots$} [l] at 103 6
\put {$\dots$} [l] at 103 48

\put {$H_{-3}$} [l] at -13 64
\put {$H_{-2}$} [l] at -13 32
\put {$H_{-1}$} [l] at -13 16
\put {$H_0$} [l] at -13 8
\put {$H_1$} [l] at -13 4
\put {$\bd^* \T$} [l] at -13 -8
\put {$\vdots$} at -10 -3
\put {$\vdots$} [B] at -10 70
\put {$\circ$} at 8 8
\put {$\omega$} at 82 82

\put {\scriptsize $0$} at 3.6 6.2
\put {\scriptsize $1$} at 12.2 6.2
\put {\scriptsize $0$} at 19.8 6.2
\put {\scriptsize $1$} at 28.4 6.2
\put {\scriptsize $0$} at 36   6.2
\put {\scriptsize $1$} at 44.2 6.2
\put {\scriptsize $0$} at 51.8 6.2
\put {\scriptsize $1$} at 60   6.2
\put {\scriptsize $0$} at 67.6 6.2
\put {\scriptsize $1$} at 76   6.2
\put {\scriptsize $0$} at 83.8 6.2
\put {\scriptsize $1$} at 92.2 6.2

\put {\scriptsize $0$} at 9 12
\put {\scriptsize $1$} at 22.5 12
\put {\scriptsize $0$} at 41 12
\put {\scriptsize $1$} at 54.5 12
\put {\scriptsize $0$} at 73 12
\put {\scriptsize $1$} at 86.5 12

\put {\scriptsize $0$} at 21 24
\put {\scriptsize $1$} at 43 24
\put {\scriptsize $0$} at 85 24

\put {\scriptsize $0$} at 45 48
\put {\scriptsize $1$} at 83 48

\put {\scriptsize $0$} at 72 75

\endpicture
$$

\vspace{.1cm}

\begin{center}
\centerline\emph{Figure 1}
\end{center}

\vspace{.1cm}
We label each edge of $\T$ by an element of $\Z_q$ such that for each vertex,
the ``downward'' edges to its $q$ successors carry labels $0, \dots, q-1$
from left to right (say), see Figure 1.  
Thus, for each $x \in  \T$, the sequence $\bigl( \sigma(n) \bigr)_{n \le 0}$
of labels on the geodesic  $\geo{x\,\om}$ has 
finite support $\{ n : \sigma(n) \ne 0 \}$. 
We write $\Sigma_q$ for the set of all those sequences.
On every horocycle, there is exactly one vertex corresponding to each 
$\sigma \in \Sigma_q$.
Thus, $\T_q$ is in one-to-one correspondence with the set 
$\Sigma_q \times \Z$, and the $k$-th horocycle is 
$H_k = \Sigma_q \times \{k\}$.

\smallskip

Now consider two trees $\T^1=\T_q$ and $\T^2=\T_r$ with roots $o_1$ 
and $o_2$ and reference ends $\om_1$ and $\om_2$, respectively.

\begin{dfn}\label{DLdef}
The Diestel-Leader graph $\DL(q,r)$ is
$$
\DL(q,r) = \{ x_1x_2 \in \T_q \times \T_r : \hor(x_1)+\hor(x_2) = 0 \}\,,
$$
and neighbourhood is given by
$$
x_1x_2 \sim y_1y_2 \iff x_1 \sim y_1 \AND x_2 \sim y_2\,.
$$
\end{dfn}

To visualize $\DL(q,r)$, draw $\T_q$ in horocyclic layers with
$\om_1$ at the top and $\bd^*\T_q$ at the bottom, and right to it $\T_r$
in the same way, but upside down, with the respective horocycles $H_k(\T_q)$
and $H_{-k}(\T_r)$ aligned. Connect the two origins $o_1$, $o_2$ by
an elastic spring. It is allowed to move along each of the two trees, 
may expand infinitely, but must always remain in horizontal position. 
The vertex set of $\DL_{q,r}$ consists of all admissible positions of 
the spring. From a position
$x_1x_2$ with $\hor(x_1) + \hor(x_2) =0$ the spring may
move downwards to one of the $r$ successors of $x_2$ in $\T_r$, and at the same 
time to the predecessor of $x_1$ in $\T_q$, or it may move upwards 
in the analogous way. Such a move corresponds
to going to a neighbour of $x_1x_2$. We see that $\DL(q,r)$ is regular
with degree $q+r$.
As the reference point in $\DL(q,r)$, we choose $o=o_1o_2$.
Figure 2 illustrates $\DL(2,2)$.

The position of $y=y_1y_2$ with respect to $x_1x_2 \in \DL(q,r)$
is described by the four numbers $\up(x_1,y_1), \dn(x_1,y_1),  \up(x_2,y_2), 
\dn(x_2,y_2)$, see below in \S \ref{renewal}, \eqref{distance} and Figure 3. 
The random walks that we are going to deal with are all such that the
transition probabilities $p(x_1x_2,y_2y_2)$ depend only on those four
parameters -- a crucial prerquisite for our approach.

$$
\beginpicture
\setcoordinatesystem units <3mm,3.5mm> 

\setplotarea x from -4 to 30, y from -3.8 to 6.4
\arrow <5pt> [.2,.67] from 4 4 to 1 7
\put{$\omega_1$} [rb] at 1.2 7.2

\put{$o_1$} [lb] at  8.15 0.2

\plot -4 -4       4 4         /         
\plot 4 4         12 -4          /      
\plot -2 -2       -2.95 -4 /            
\plot -.5 -2      -1.9 -4     /         
\plot -.5 -2      -.85 -4   /           
\plot 1 -2        .2 -4      /          
\plot 1 -2        1.25  -4     /        
\plot 2.5 -2      2.3  -4     /         
\plot 2.5 -2      3.35 -4      /        
\plot 5.5 -2      4.65  -4    /         
\plot 5.5  -2     5.7  -4      /        
\plot 7  -2       6.75   -4   /         
\plot 7 -2        7.8  -4      /        
\plot 8.5  -2     8.85  -4    /         
\plot 8.5  -2     9.9  -4      /        
\plot 10  -2      10.95  -4   /         
\plot 0  0        -.5  -2      /        
\plot 2 0         1 -2     /            
\plot 2 0         2.5   -2     /        
\plot 6 0         5.5    -2    /        
\plot 6 0         7 -2         /        
\plot 8 0         8.5 -2       /        
\plot 2 2         2 0         /         
\plot 6 2         6 0         /         

\arrow <5pt> [.2,.67] from 22 -4 to 25 -7
\put{$\omega_2$} [lt] at 25.2 -7.2

\put{$o_2$} [rt] at  17.95 -.2

\plot 14  4       22 -4       /         
\plot 22 -4         30 4         /      
\plot 16 2       15.05 4 /              
\plot 17.5 2     16.1  4     /          
\plot 17.5 2     17.15  4   /           
\plot 19  2       18.2  4      /        
\plot 19 2         19.25   4     /      
\plot 20.5 2       20.3   4     /       
\plot 20.5  2      21.35 4      /       
\plot 23.5 2       22.65  4    /        
\plot 23.5  2      23.7   4      /      
\plot 25   2       24.75   4   /        
\plot 25  2        25.8  4      /       
\plot 26.5   2     26.85  4    /        
\plot 26.5   2     27.9   4      /      
\plot 28   2      28.95   4   /         
\plot 18  0        17.5  2      /       
\plot 20 0         19  2     /          
\plot 20 0         20.5   2     /       
\plot 24 0         23.5    2    /       
\plot 24 0         25 2         /       
\plot 26 0         26.5 2       /       
\plot 20 -2        20 0         /       
\plot 24 -2        24 0         /       
\put {$\circ$} at 8 0
\put {$\circ$} at 18 0
\plot 8.25 0  12.1 0 /
\plot 13.9 0 17.78 0 /
\plot    12.1   0    12.25 .4    12.25 -.4   12.55 .4   12.55 -.4
        12.85 .4   12.85 -.4  13.15 .4   13.15 -.4  13.45 .4
        13.45 -.4  13.75 .4   13.75 -.4  13.9 0     13.9  0 /

\setdashes <2pt>
\putrule from -4.5 -7  to  12.5 -7
\putrule from  13.5 7  to  30.5 7

\put {$\bd^*\T^1$} [r] at -5 -7
\put {$\bd^*\T^2$} [l] at 31 7

\put {$\vdots$} at 4 -5.2
\put {$\vdots$} at 22 5.5

\endpicture
$$

\begin{center}
\emph{Figure 2}
\end{center}

\vspace{.4cm}

We now recall in more detail the construction of the \emph{lamplighter group} 
$\Z_q \wr \Z$. The group of all finitely supported configurations is
$$
\CC = \{ \eta: \Z \to \Z_q \,,\;\supp(\eta) \;\ \mbox{finite}\, \}
$$
with pointwise addition modulo $q$. The group $\Z$ acts on $\CC$
by translations $k \mapsto T_k:\CC \to \CC$ with $T_k\eta(m) = \eta(m-k)$.
The resulting semidirect product $\Z \rightthreetimes \CC$ is
$$
\Z_q \wr \Z = \{ (\eta,k) : \eta \in \CC\,,\;k \in \Z\} 
\quad \mbox{with group operation}\quad
(\eta,k)(\eta',k') = (\eta + T_k\eta', k+k')
$$
We identify each $(\eta,k) \in \Z_q \wr \Z$ with the vertex 
$x_1x_2 \in \DL(q,q)$, where according to the identification 
$\T_q \lra \Sigma_q \times \Z$,
the vertices $x_i$  are given by 
\begin{equation}\label{identif}
\begin{gathered}
x_1 = (\eta_k^-,k) \AND x_2=(\eta_k^+,-k)\,,\quad\mbox{where}\\
\eta_k^- = \eta|_{(-\infty\,,\,k]} \AND 
\eta_k^+ = \eta|_{[k+1\,,\,\infty)}\,,
\end{gathered}
\end{equation}
both written as sequences over the non-positive integers.

This is a one-to-one correspondence between $\Z_q \wr \Z$ and 
$\DL(q,q)$, and that group acts transitively and fixed-point-freely 
on the graph. Namely, the action of $m \in \Z$ is given by 
$x_1x_2 = (\sigma_1,k)(\sigma_2,-k) 
\mapsto y_1y_2 = (\sigma_1,k+m)(\sigma_2,-k+m)$,
and the action of the group of configurations is pointwise addition modulo $q$. 
Write $\de_k^l$ for the configuration in $\CC$ with value 
$l$ at $k$ and $0$ elsewhere.
Then $\DL(q,q)$ is the (right) Cayley graph of $\Z_q \wr \Z$ with respect to
the symmetric set of generators
$$
\{ (\de_1^l,1)\,,\; (\de_0^l,-1) : l \in \Z_q \}\,,
$$
i.e., an edge corresponds to multiplying with a generator on the right.
(This is precisely the set of generators considered in \cite{GriZuk}
and \cite{DiSc} when computing the spectrum of the associated 
SRW-operator.)

\smallskip

Returning to $\DL=\DL(q,r)$, the transition matrix $P_{\al}$ of the
random walk that we have described in the Introduction is given as follows.
For $x=x_1x_2, y=y_1y_2 \in \DL(q,r)$
\begin{equation}\label{random-walk}
p_{\al}(x,y) = \begin{cases} \al/q 
                   & \text{if}\; y_1^- = x_1 \;\text{and}\;y_2=x_2^-\\
                                 (1-\al)/r 
                   & \text{if}\; y_1 = x_1^- \;\text{and}\;y_2^-=x_2\\
                 0 & \text{otherwise.}
                   \end{cases}
\end{equation}

\section{Simple random walk with drift on a homogeneous tree}\label{trees}

In general, if $P$ is a transition matrix over a set $X$ and $\{X_i : i \in I\}$
is a partition of $X$ with the associated quotient map $\pi: X \to I$,
then one says that $P$ \emph{factorizes} (or \emph{projects}) with respect 
to $\pi$, if $\;\wt p(i,j) := \sum_{y \in X_j} p(x,y)\;$ does not depend on the 
specific choice of $x \in X_i$. In this case, the Green kernel $\wt G$ 
associated with $\wt P = \pi(P)$ also satisfies
\begin{equation}\label{green-project}
\wt G(i,j) = \sum_{y \in X_j} G(x,y)\,,\quad x \in X_i\,.
\end{equation}
In our case, we have three natural, neighbourhood preserving projections
$\pi_1: \DL \to \T_q\,$, $\pi_2: \DL \to \T_r\,$, and $\wt\pi: \DL \to \Z\,$,
given by
$$\pi_1(x_1x_2)=x_1\,,\quad \pi_2(x_1x_2)=x_2\,, \AND
\wt\pi(x_1x_2)=\hor(x_1)\,.
$$
$P_{\al}$ factorizes with respect to each of them. Let 
$\pi_1(P_{\al}) = P_1$,  $\pi_2(P_{\al}) = P_2$ and
$\wt\pi(P_{\al}) = \wt P$. Then  $P_1=P_{\al,q}$ on $\T^1 = \T_q$, 
$P_2 = P_{1-\al,r}$ on $\T^2 = \T_r$, and $\wt P = \wt P_{\al}$ on $\Z$,
where
\begin{equation}\label{projections}
p_{\al,q}(x_1,y_1) = \begin{cases} \al/q & \text{if}\; y_1^- = x_1 \\
                             (1-\al) & \text{if}\; y_1 = x_1^- \\
			     0 & \text{otherwise,}
                                    \end{cases}
\qquad
\wt p_{\al}(k,l) = \begin{cases} \al & \text{if}\; l=k+1 \\
                             (1-\al) & \text{if}\; l=k-1\\
			     0 & \text{otherwise.}
                                    \end{cases}
\end{equation}
The projected random walks are well understood. Everybody is familiar with the
gambler's process $\wt P_{\al}$ on $\Z$. We outline
the basic features of $P_{\al,q}$ on $\T_q$ (or, equivalently,
$P_{1-\al,r}$ on $\T_r$).

\smallskip

{\bf Spectral radius.} The spectral radius of any irreducible transition matrix
is defined as $\rho(P) = \limsup_n p^{(n)}(x,y)^{1/n}$. It is independent
of $x, y$. In the specific case of our random walks with drift parameter
$\al$, we have
\begin{equation}\label{radii}
\rho(P_{\al})_{\DL} = \rho(P_{\al,q})_{\T_q} = 
\rho(P_{1-\al,r})_{\T_r} = \rho(\wt P_{\al})_{\Z} = 2\sqrt{\al(1-\al)}\,.
\end{equation}
(The subscript refers to the respective underlying graph.)
For $\wt P_{\al}$ on $\Z$, this is well known. For $P_{\al,q}$ on $\T_q$, 
it can be easily computed in various ways. See e.g. 
{\sc Saloff-Coste and Woess} \cite{SCW}, Example 1.

\smallskip

{\bf Green kernel.} The -- simple -- computations of the Green kernel
$G_1 = G_{\al,q}$ associated with $P_{\al,q}$ can be done following the method 
of \S 1.D in \cite{Wbook}, see also \cite{Woe}, (3.9). The main point ist
that we have a nearest neighbour random walk on a tree (transition
probabilities are positive only between neighbours). Thus, if $F_1(x_1,y_1)$ is 
the probability that the random walk starting at $x_1$ ever hits $y_1$ 
($x_1, y_1 \in \T_q$), then
\begin{equation}\label{tree-lemma}
F_1(x_1,y_1) = F_1(x_1,w_1)F_1(w_1,y_1) \quad \mbox{for all}\; 
w_1 \in \geo{x_1\,y_1}\,.
\end{equation}
Furthermore, since $p_1(x_1,y_1)$ depends only on $\up(x_1,y_1)$ and 
$\dn(x_1,y_1)$, the same is true for $F_1(x_1,y_1)$ and $G_1(x_1,y_1)$. 
In particular,
$$
F_1^- = F_1(x_1,x_1^-) \AND F_1^+ = F_1(x_1^-,x_1)
$$
are independent of $x_1 \in \T_q$ as well as $G_1(x_1,x_1)$.
Using these facts, and setting $\al^+=\max\{\al, 1-\al\}$, one computes
\begin{equation}\label{G1xy}
\begin{aligned}
G_1(x_1,y_1) &= F_1(x_1,y_1)\,G_1(y_1,y_1) 
= (F_1^{-})^{\up(x_1,y_1)}(F_1^+)^{\dn(x_1,y_1)}
\frac{q}{(q+1)\al^+ -1}\,, \quad\mbox{where}\\[7pt]
F_1^- &= \begin{cases} \dps \frac{1-\al}{\al} 
                           &\dps \text{if}\; \al \ge \frac12\,,   \\[7pt]
                           \dps 1& \dps \text{if}\; \al \le \frac12\,, 
                                    \end{cases} 
\qquad
F_1^+ = \begin{cases} \dps \frac{1}{q} 
                           & \dps \text{if}\; \al \ge \frac12\,,    \\[7pt]
                           \dps   \frac{\al}{(1-\al)q} 
			   & \dps \text{if}\; \al \le \frac12\,.
                                    \end{cases}\\
\end{aligned}
\end{equation}

{\bf Martin compactification.} By \eqref{tree-lemma}, the Martin kernel 
$K_1 = K_{\al,q}$ associated with $P_{\al,q}$ satisfies
$$
K_1(x_1,y_1) = \frac{F_1(x_1,y_1)}{F_1(o_1,y_1)} 
= \frac{F_1(x_1,c_1)}{F_1(o_1,c_1)}\,,\quad \mbox{where}\; 
c_1 = x_1 \wedge y_1
$$
(the confluent with respect to $o_1$). From here, the following
is almost immediate.

\begin{pro}\label{tree-martin}
The Martin compactification of\/ $\T_q$ with respect to $P_{\al,q}$
is the end compactification\/ $\wh \T_q$, and for $\xi_1 \in \bd \T_q$, the 
Martin kernel is given by $K_1(x_1,\xi_1) = K_1(x_1,c_1)\,$, 
where $c_1 = x_1 \wedge \xi_1\,$.

Furthermore, each function $K_1(\cdot,\xi_1)$, $\xi_1 \in \bd \T_q$, 
is minimal harmonic for $P_{\al,q}$.
\end{pro}

For general transient nearest neighbour random walks on arbitrary locally
finite trees, this is due to by {\sc Cartier} \cite{Car}, and in the
specific case of free groups (which is close to, but not identical with
our situation), it was shown previously by 
{\sc Dynkin and Malyutov}~\cite{DyMa}.

The analogous results for $P_{1-\al,r}$ on $\T_r$ are obtained
from the above by exchanging $\al$ with $1-\al$ and $q$ with $r$.
When $\al \ne 1/2$, the same computations are also valid for 
$\wt P_{\al}$ on $\Z$, setting $q=1$. When $\al = 1/2$ then
$\wt P_{\al}$ is of course \emph{recurrent}, i.e., the
associated Green kernel diverges.

Below in \S \ref{martin}, we shall also need the following 
functions on $\T_q \times \T_q$, which we call \emph{(generalized) 
spherical functions.} We set $\up = \up(x_1,y_1)$, $\dn = \dn(x_1,y_1)$ 
and $\be = (1-\al)/\al$.
\begin{equation}\label{spherical}
\phi_{\al,q}(x_1,y_1) = \begin{cases}
\dfrac{1}{q^{\dn}}\left( \dfrac{q+1}{q-1} + d(x_1,y_1)\right)\,,
&\mbox{if }\; \al=\dfrac12\,,\\[10pt]
\dfrac{1}{(q\,\be^2)^{\dn}}\left( \dfrac{q\,\be+1}{q\,\be^2-1} + 
\dfrac{\be^{\up}-1}{\be -1}+ \dfrac{\be^{\dn}-1}{\be -1}
\right)&\mbox{if }\; \al < \frac12\,,\\[10pt]
\dfrac{\be^{d(x_1,y_1)}}{q^{\dn}}
\left( \dfrac{q\,\be^{-1}+1}{q\,\be^{-2}-1} + 
\dfrac{\be^{-\up}-1}{\be^{-1} -1}+ 
\dfrac{\be^{-\dn}-1}{\be^{-1} -1}
\right)&\mbox{if }\; \al > \frac12\,.
\end{cases}
\end{equation}
(Recall that $d(x_1,y_1) = \up+\dn\,$.) Then $\phi_{\al,q}(\cdot,y_1)$ is 
$P_{\al,q}$-harmonic on $\T_q$ for each $y_1 \in \T_q$.

\section{Green kernel asymptotics}\label{renewal}

We now embark on the main computational part of this paper. 
We consider $P_{\al}$ on $\DL = \DL(q,r)$, and 
\emph{\sc we shall always assume that $\al \le 1/2$}, since all results 
in the case $\al \ge 1/2$ are obtained 
from the former by exchanging the role of the two trees (i.e., exchanging $r$
with $q$). 

We want to derive asymptotic estimates of the associated Green kernel 
$G(x,y) = G_{\al}(x,y)$, where $x=x_1x_2$ and $y=y_1y_2 \in DL$
and the graph distance $d(x,y) \to \infty$. The latter means that 
at least one of $d(x_1,y_1)$ and $d(x_2,y_2)$ (distances in the
respective trees) tends to $\infty$. 
We remark here that 
\begin{equation}\label{distance}
\begin{aligned}
d(x,y) &= d(x_1,y_1) + d(x_2,y_2) - |\hor(y_1)-\hor(x_1)|\,,\\
d(x_i,y_i) &= \up_i+\dn_i\,, \quad
\hor(y_i)-\hor(x_i) = \dn_i - \up_i\,,\quad(i=1,2), \AND\\
\up_1 + \up_2 &= \dn_1 + \dn_2\,,
\quad\mbox{where}\quad \up_i= \up(x_i,y_i) \AND \dn_i = \dn(x_i,y_i)\,.
\end{aligned}
\end{equation}
(Cf. {\sc Bertacchi} \cite{Ber} for the distance formula.) 
In terms of the lamplighter moving along $\Z$ (with the lamps -- possibly red 
\emph{and} green -- sitting on the edges, as described in the Introduction),
$\up_1$ is the minimal number of steps the lamplighter has to walk
in the negative direction in order to obtain the new position and configuration
encoded in the vertex $y=y_1y_2$ of $\DL$, and $\up_2$ is analogous in
the positive direction. 

We set $c_i = x_i \cf y_i$.
See Figure 3. We also choose $a_i, b_i \in \T^i$ with $x_i \lle a_i$,
$y_i \lle b_i$, such that $\hor(a_1) = \hor(b_1)=-\hor(c_2)$
and $\hor(a_2) = \hor(b_2)=-\hor(c_1)$, i.e., the pairs
$a_1c_2$, $b_1c_2$, $c_1a_2$, $c_1b_2$ belong to $\DL$.
In particular, $d(x,y) \to \infty$ means that $\spn \to \infty$, where
$$
\spn = \spn(x,y) = \up_1+\up_2=\dn_1+\dn_2= -\hor(c_1)-\hor(c_2)
$$
is the \emph{span} of $x$ and $y$.

$$
\beginpicture
\setcoordinatesystem units <3.2mm,3.7mm> 

\setplotarea x from -4 to 24, y from -3.8 to 6.4

\put{$x_1$} [lb] at  8.15 0.2
\put{$y_1$} [rb] at 1.85 2.2
\put{$c_1$} [lb] at 4.15 4.2
\put{$\omega_1$} [rb] at 0.2 8.2 
\put{$a_1$} [lt] at 12 -4.2
\put{$b_1$} [rt] at -3.8 -4.2

\arrow <5pt> [.2,.67] from 4 4 to 0 8 
\plot 8 0  4 4  2 2 /

\put{$x_2$} [rt] at  11.95 -.2
\put{$y_2$} [lt] at  22.2  1.8
\put{$c_2$} [rt] at  15.95 -4.2
\put{$\omega_2$} [lt] at 20.2 -8.2
\put{$a_2$} [rb] at  7.8 4.2
\put{$b_2$} [lb] at  24 4.2


\plot 12 0   16 -4   22 2 /         
\arrow <5pt> [.2,.67] from 16 -4 to 20 -8

\multiput {$\scs\bullet$} at 
                            4 4  16 -4 
                            -4 -4  12 -4  8 4  24 4 /
\multiput {$\bullet$} at 8 0  12 0  2 2  22 2  /

\setdashes <3pt>

\plot  2 2   -4 -4 /   
\plot  8 0   12 -4 /
\plot  12 0   8 4 /
\plot  22 2  24 4 /


\setdots<4pt>
\plot -9.5 0  12 0 /
\plot -6.5 2  22 2 /
\plot -9.5 4  24 4 /
\plot -9.5 -4  16 -4 /

\put{$\up_1\!\!\left\{\rule[-6mm]{0mm}{0mm}\right.$}[rb] at -9.5 -.1 
\put{$\up_2\!\!\left\{\rule[9mm]{0mm}{0mm}\right.$}[rt] at -9.5 0.1
\put{$\dn_1\!\!\left\{\rule[5mm]{0mm}{0mm}\right.$}[rb] at -6.5 2
\put{$\dn_2\!\!\left\{\rule[-9.95mm]{0mm}{0mm}\right.$}[rt] at -6.5 2.1

\endpicture
$$
\vspace{.1cm}

\begin{center}
\emph{Figure 3}
\end{center}

\vspace{.4cm}

The following is the first main result of this paper.

\begin{thm}\label{estimates} Referring to \eqref{distance} and Figure 3,
suppose that $d(x,y) \to \infty$, and hence $\spn = \spn(x,y) \to \infty$.
Then we have the following.\\[3pt]
{\rm (a)} If $\al > 1/2$ and $\be = (1-\al)/\al$ then
$$
\begin{aligned}
G(x,y) \sim \frac{A_{\be}}{(q\,\be)^{\dn_1}r^{\dn_2}}&\left(
B_{\be}\, \frac{\be^{\spn}-\be^{\up_1}}{\be^{\spn}-1}\,
        \frac{\be^{\spn}-\be^{\dn_1}}{\be^{\spn}-1}\,\frac{1}{\be^{\spn}}
+   	\frac{\be^{\spn}-\be^{\up_1}}{\be^{\spn}-1}\,
        \frac{\be^{\dn_1}-1}{\be^{\spn}-1}\right.\\
&\quad+   	\left.\frac{\be^{\up_1}-1}{\be^{\spn}-1}\,
        \frac{\be^{\spn}-\be^{\dn_1}}{\be^{\spn}-1}\,
+B_{\be}^*\, \frac{\be^{\up_1}-1}{\be^{\spn}-1}\,
        \frac{\be^{\dn_1}-1}{\be^{\spn}-1} \right),
\end{aligned}
$$
where
$$
A_{\be} = \frac{G_1(o_1,o_1)\,G_2(o_2,o_2)}{\wt G(0,0)} 
= \frac{q\,r\,(\be^2-1)}{(q\,\be -1)(q\,r-1)}\,.
$$\\
{\rm (b)} If $\al = 1/2$ then
$$
G(x,y) \sim  \,\frac{A_1}{\spn^4\,q^{\dn_1}\,r^{\dn_2}}
\left(\frac{q+1}{q-1}\,\up_2\,\dn_2 +\spn\,\up_2\,\dn_1 + 
\spn\,\up_1\,\dn_2 + \frac{r+1}{r-1}\,\up_1\,\dn_1 \right),
$$
where 
$$
A_1 = \frac{G_1(o_1,o_1)\,G_2(o_2,o_2)}{2} = \frac{2\,q\,r}{(q-1)(r-1)}\,.
$$
\end{thm} 

According to the way how $y$ tends to infinity geometrically (when we
think of $x$ being fixed), one or more of the four terms will
dominate the others, as we shall see below. 

As mentioned at the beginning, the case $\al < 1/2$ is obtained by 
exchanging $r \leftrightarrow q$ and $\al \lra 1-\al$. Equivalently,
we may use Lemma \ref{exchange} and apply statement (a) 
of Theorem \ref{estimates} to $G^*(x,y)$, with $\be^* = 1/\be$ in
the place of $\be$. 

We now start make our (laborious) way towards the proof of 
Theorem \ref{estimates}. The following is obvious, but crucial.

\begin{lem}\label{coordinates} The Green kernel $G(x,y)$ depends only
on $\up_1, \dn_1, \up_2, \dn_2$.
\end{lem}

Let $Z_n$ be the random position of the $P_{\al}$-walk. This is a $\DL$-valued
random variable defined on a suitable probability space (trajectory space).
We write $\Prb_x = \Prb[\,\,\cdot \mid Z_0=x]$ and $\Ex_x$ for the associated 
expectation. Also, $\1_{[\,\cdots]}$ will denote the indicator
function of an event $[\,\cdots]$ in the trajectory space. 
The projection $Z_n^i= \pi_i(Z_n)$ represents the random position at time
$n$ of the $P_i$-walk on $\T^i$, $i=1,2$, and the random variable
$\wt Z_n = \wt\pi(Z_n)$ realizes the $n$-th position of the
$\wt P_{\al}$-walk on $\Z$.

We shall use several \emph{stopping times.} If $x=x_1x_2 \in \DL$,
$x_i \in \T^i$ ($i=1,2$), resp. $k \in Z$, then we set
$$
\begin{gathered}
\tm(x) = \inf \{ n \ge 0 : Z_n = x \}\,,\quad
\tm_i(x_i)  = \inf \{ n \ge 0 : Z_n^i = x_i \}\quad(i=1,2)\,,\AND\\
\wt \tm(k) = \inf \{ n \ge 0 : \wt Z_n = k \}\,.
\end{gathered}
$$
Note that these random variables are all defined on the same trajectory
space associated with $P_{\al}$.
\begin{lem}\label{stop}
Referring to the situation of Figure 3, we have
$$
\tm_1(c_1) = \wt\tm(-\up_1) \AND \tm_2(c_2) = \wt\tm(\up_2)\quad 
\Prb_x\mbox{-almost surely.}
$$
Furthermore, in order to reach $y$ starting in $x$, both $\Z_n^i$
have to pass through $c_i$, $i=1,2$, i.e.,
$$
\max \{ \tm_1(c_1), \tm_2(c_2) \} \le \tm(y) \quad 
\Prb_x\mbox{-almost surely.}
$$
\end{lem}

\begin{proof} The $P_{\al}$-walk on $\DL$ as well as the projected
random walks are of nearest neighbour type. Thus, starting in $x$, 
the first point in the set $\{ v=v_1v_2 \in \DL : \hor(v_1)=\hor(c_1) \}$
visited by $Z_n$ must be of the form $c_1v_2$. This translates
to $\tm_1(c_1) = \wt\tm(\up_1)$, and exchanging the roles of the
two trees, also to $\tm_2(c_2) = \wt\tm(\up_2)$.
The same type of argument shows the last statement.  
\end{proof}

The last lemma leads us to the identities
\begin{equation}\label{c1c2}
\Prb_x[ \tm_1(c_1) < \tm_2(c_2) ] = \varphi_1(\up_1,\up_2)
\AND
\Prb_x[ \tm_2(c_2) < \tm_1(c_1) ] = \varphi_2(\up_1,\up_2)\,,
\end{equation}
where for $k, l \ge 0$, the probability that the $\wt P_{\al}$-walk
on $\Z$ starting in $0$ reaches $-k$ before $l$ is $\varphi_1(k,l)$,
and the probability that it reaches $l$ before $-k$ is
$\varphi_2(k,l) = 1-\varphi_1(k,l)$. It is a well-known exercise to
compute these functions, and they are given by
\begin{equation}\label{varphi}
\begin{alignedat}{3}
\varphi_1(k,l) &= \frac{\be^{k+l} - \be^k}{\be^{k+l} - 1}&\AND
\varphi_2(k,l) &= \frac{\be^{k} - 1}{\be^{k+l} - 1}\,,
&\quad&\mbox{with}\;\; \be=\frac{1-\al}{\al}\,,\quad
                      \mbox{if}\;\;\al\ne \frac12\,;\\[6pt]
\varphi_1(k,l) &=\frac{l}{k+l}&\AND		      
\varphi_2(k,l) &=\frac{k}{k+l}\,,&\quad&\mbox{if}\;\;\al = \frac12\,.		      
\end{alignedat}
\end{equation} 
See e.g. {\sc Kemeny and Snell} \cite{KeSn}, \S 7.1, in particular (5) and (6)
on p. 153.
Next, let us introduce the function
\begin{equation}\label{psi}
\psi(k) = \left( \frac{\al\, r}{(1-\al)\,q}\right)^{\! k}\,,\quad k \in \Z\,.
\end{equation}
If we set $\mm(x) = \psi\bigl(\hor(x_1)\bigr)$, where $x=x_1x_2 \in \DL$,
then we have $\mm(x)\, p_{\al}(x,y) = \mm(y) \,p_{\al}(y,x)$ for all
$x, y \in \DL$. That is, $P_{\al}$ is \emph{$\mm$-reversible}, and
we also get
\begin{equation}\label{reverseG}
G(x,y) = \psi\bigl( \hor(y_1) - \hor(x_1) \bigr) \,G(y,x) \quad
\mbox{for all}\;\; x=x_1x_2\,,\ y=y_1y_2 \in \DL\,.
\end{equation}

\begin{pro}\label{decomp} Referring to the situation of Figure 3, 
we have the following decomposition.
$$
\begin{alignedat}{3}
& \hspace{3cm} & G(x,y) =\;&\varphi_1(\up_1,\up_2) \,\varphi_1(\dn_1,\dn_2)\,\psi(\dn_1)\,
          G(c_1b_2,c_1a_2)& \hspace{3.2cm} \mbox{\rm (I)}& \\
& &        +\ &\varphi_1(\up_1,\up_2) \,\varphi_2(\dn_1,\dn_2)\,\psi(\dn_1)\,
        G(b_1c_2,c_1a_2)& \mbox{\rm (II)}& \\ 
& &        +\ &\varphi_2(\up_1,\up_2) \,\varphi_1(\dn_1,\dn_2)\,\psi(-\dn_2)\,
          G(c_1b_2,a_1c_2)& \mbox{\rm (III)}& \\ 
& &       +\ &\varphi_2(\up_1,\up_2) \,\varphi_2(\dn_1,\dn_2)\,\psi(-\dn_2)\,
          G(b_1c_2,a_1c_2)\,.& \mbox{\rm (IV)}& 
\end{alignedat}
$$
\end{pro}	   

\begin{proof} By \eqref{stop} and \eqref{c1c2}, we have
$$
G(x,y) 
= \Ex_x\Bigl( \1_{[ \tm_1(c_1) < \tm_2(c_2) ]}\,G(Z_{\tm_1(c_1)},y)\Bigr)  
+ \Ex_x\Bigl( \1_{[ \tm_2(c_2) < \tm_1(c_1) ]}\,G(Z_{\tm_2(c_2)},y)\Bigr)\,.
$$
If $\tm_1(c_1) < \tm_2(c_2)$ and $Z_{\tm_1(c_1)} =c_1w_2$, then we must have
$c_2 \lle w_2$ and futhermore $w_2 \cf b_2 = c_2$. Thus,
$\up(w_2,y_2) = \up(a_2,y_2)$ and $\dn(w_2,y_2) = \dn(a_2,y_2)$.
Lemma \ref{coordinates} implies $G(Z_{\tm_1(c_1)},y) = G(c_1a_2,y)$.
In the same way, $G(Z_{\tm_2(c_2)},y) = G(a_1c_2,y)$. Thus,
$$
\begin{aligned}
G(x,y) 
&= \Prb_x[ \tm_1(c_1) < \tm_2(c_2) ]\,G(c_1a_2,y) +
+ \Prb_x[ \tm_2(c_2) < \tm_1(c_1) ]\,G(a_1c_2,y)\\
&= \varphi_1(\up_1,\up_2)\,G(c_1a_2,y)+\varphi_2(\up_1,\up_2)\,G(a_1c_2,y)\,.
\end{aligned}
$$
Using \eqref{reverseG}, we get $G(c_1a_2,y) = \psi(\dn_1)\,G(y,c_1a_2)$.
Applying once more \eqref{stop} and \eqref{c1c2},
$$
G(y,c_1a_2)
=\Ex_y\Bigl( \1_{[ \tm_1(c_1) < \tm_2(c_2) ]}\,G(Z_{\tm_1(c_1)},y)\Bigr)  
+ \Ex_y\Bigl( \1_{[ \tm_2(c_2) < \tm_1(c_1) ]}\,G(Z_{\tm_2(c_2)},y)\Bigr)\,.
$$
We can repeat the above argument with $y$ in the place of $x$ and
$c_1a_2$ in the place of $y$, and we have to replace $a_1, a_2$ with 
$b_1, b_2$. Therefore
$$
G(y,c_1a_2)
= \varphi_1(\dn_1,\dn_2)\,G(c_1b_2, c_1a_2) +
\varphi_2(\dn_1,\dn_2)\,G(b_1c_2, c_1a_2)\,.
$$
Analogously, $G(a_1c_2,y) = \psi(-\dn_2)\,G(y,a_1c_2)$ and
$$
G(y,a_1c_2)
= \varphi_1(\dn_1,\dn_2)\,G(c_1b_2, a_1c_2) +
\varphi_2(\dn_1,\dn_2)\,G(b_1c_2, a_1c_2)\,. 
$$
Combining these formulas, we obtain the proposed decomposition.
\end{proof}

Thus, in order to understand the asymptotics of $G(x,y)$ in the
general case of Figure 3, we can reduce our computations to the following four
basic cases of relative positions of $x$ and $y$.
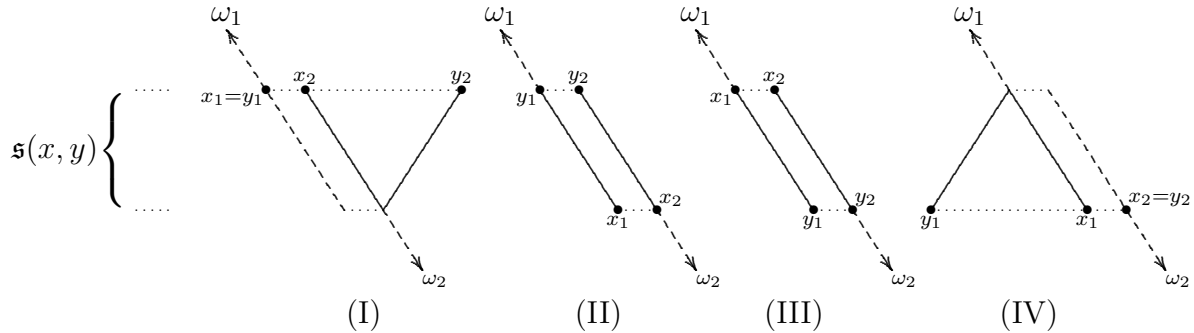
$$
\beginpicture
\setcoordinatesystem units <2.6mm,4mm> 

\setplotarea x from -10 to 41, y from -3 to 6

\plot -2 4  2 0  6 4 /

\plot 10 4  14 0 /
\plot 12 4  16 0 /

\plot 20 4  24 0 /
\plot 22 4  26 0 /

\plot 30 0  34 4  38 0 /

\arrow <5pt> [.2,.67] from -5.6 5.6 to -6 6
\arrow <5pt> [.2,.67] from 3.6 -1.6 to 4 -2

\arrow <5pt> [.2,.67] from  8.4 5.6 to 8 6
\arrow <5pt> [.2,.67] from 17.6 -1.6 to 18 -2

\arrow <5pt> [.2,.67] from  18.4 5.6 to 18 6
\arrow <5pt> [.2,.67] from 27.6 -1.6 to 28 -2

\arrow <5pt> [.2,.67] from  32.4 5.6 to 32 6
\arrow <5pt> [.2,.67] from  41.6 -1.6 to 42 -2

\setdashes <2.5pt>
\plot 0 0  -5.6 5.6 /
\plot 2 0  3.6 -1.6 /

\plot 10 4  8.4 5.6 /
\plot 16 0  17.6 -1.6 /

\plot 20 4  18.4 5.6 /
\plot 26 0  27.6 -1.6 /

\plot 34 4  32.4 5.6 /
\plot 36 4  41.6 -1.6 /

\setdots<3pt>
\plot -4 4  6 4 /
\plot  0 0  2 0 /

\plot 10 4  12 4 /
\plot 14 0  16 0 /

\plot 20 4  22 4 /
\plot 24 0  26 0 /

\plot 34 4  36 4 /
\plot 30 0  40 0 /

\plot -9 4  -11 4 /
\plot -9 0  -11 0 /

\put{$\scs x_1=y_1$} [rt] at  -4.1 3.9
\put{$\scs x_2$} [b] at -2 4.2
\put{$\scs y_2$} [b] at 6 4.2

\put{$\scs x_1$} [t] at 14 -.2
\put{$\scs y_1$} [rt] at 9.9 3.9
\put{$\scs x_2$} [lb] at 16.1 0.15
\put{$\scs y_2$} [b] at 12 4.2

\put{$\scs y_1$} [t] at 24 -.2
\put{$\scs x_1$} [rt] at 19.9 3.9
\put{$\scs y_2$} [lb] at 26.1 0.15
\put{$\scs x_2$} [b] at 22 4.2

\put{$\scs x_1$} [t] at 38 -.2
\put{$\scs y_1$} [t] at 30  -.2
\put{$\scs x_2=y_2$} [lb] at 40.1 0.15

\multiput{$\om_1$} [b] at -6 6.15  8 6.15  18 6.15  32 6.15 /
\multiput{$\scs \om_2$} [lt] at 4 -2.15  18 -2.15  28 -2.15  42 -2.15 /

\put{(I)} at  1 -3.5
\put{(II)} at 13 -3.5
\put{(III)} at 23 -3.5
\put{(IV)} at 35 -3.5

\multiput {$\scs\bullet$} at -4 4  -2 4  6 4    
                             10 4  14 0  12 4  16 0  
			     20 4  24 0  22 4  26 0  
			     30 0  38 0  40 0 /

\put{$\spn(x,y)\!\left\{\rule[8.48mm]{0mm}{1mm}\right.$}[rt] at -11 4.1

\endpicture
$$
\vspace{.1cm}

\begin{center}
\emph{Figure 4}
\end{center}

\vspace{.4cm}
In all four cases, $\spn = \spn(x,y) \to \infty$. In case (I), $\up_1=\dn_1=0$
and $\up_2=\dn_2=\spn$. In case (II), $\up_1=\dn_2=\spn$ and $\dn_1=\up_2=0$.
In case (III), $\dn_1=\up_2=\spn$ and $\up_1=\dn_2=0$. In case (IV), 
$\up_1=\dn_1=\spn$ and $\up_2=\dn_2=0$.

We start with a extended version of case II, see Figure 5.
$$
\beginpicture
\setcoordinatesystem units <2.2mm,2.56mm> 

\setplotarea x from -4 to 24, y from -3.8 to 6.4

\put{$y_1$} [rb] at 1.85 2.2
\put{$c_1$} [lb] at 4.15 4.2
\put{$\omega_1$} [rb] at 0.7 7.7 
\put{$x_1$} [rt] at 12 -4.2

\arrow <5pt> [.2,.67] from 12 -4 to 0.5 7.5 
\plot 4 4  2 2 /

\put{$x_2$} [rt] at  15.95 -4.2
\put{$\omega_2$} [lt] at 19.7 -7.7
\put{$y_2$} [rb] at  9.8 2.2
\put{$w_2$} [lb] at  26 2.2
\put{$v_2 = v(w_2)$} [lt] at  20 -3.8

\arrow <5pt> [.2,.67] from 10 2 to 19.5 -7.5

\multiput {$\scs\bullet$} at 2 2  4 4  12 -4  10 2  16 -4  20 -4  26 2 /

\setdashes <3pt>                             
\plot 18 -6  26  2 /

\setdots<4pt>

\plot -6.5 2  26 2 /
\plot -11 4  4 4 /
\plot -11 -4  20 -4 /

\put{$\spn= \up_1\!\!\left\{\rule[11mm]{0mm}{0mm}\right.$}[rt] at -11 3.95
\put{$\dn_1\!\!\left\{\rule[3.5mm]{0mm}{0mm}\right.$}[rb] at -6.5 2
\put{$\dn_2\!\!\left\{\rule[-7.34mm]{0mm}{0mm}\right.$}[rt] at -6.5 1.9

\endpicture
$$
\begin{center}
\emph{Figure 5}
\end{center}

\begin{pro}\label{IIbis} If, as in Figure 5, $x=x_1x_2$ and $y=y_1y_2$ 
with $\spn = \spn(x,y)$ are such that $\up_1 - \dn_1 = \dn_2 \to \infty$, 
$\up_2 = 0$ and $\dn_1$ is arbitrary, then
\begin{equation}\label{vertical}
G_1(x_1,y_1) = C(\spn)  \, r^{\dn_2}\, G(x,y) + R(\dn_1,\dn_2)\,,
\end{equation}
where
$$
C(\spn) \to \frac{\wt G(0,0)}{G_2(o_2,o_2)}\quad
\mbox{if}\;\;\al \ne \frac12\,,\qquad\quad
C(\spn)  \sim \frac{2\,\spn}{G_2(o_2,o_2)}\quad
\mbox{if}\;\;\al = \frac12\,,
$$
and $0 < R(\dn_1,\dn_2) < G_1(x_1,y_1)$ with 
$$
\lim_{\dn_2 \to \infty} R(\dn_1,\dn_2) = 0\,.
$$
\end{pro}

\begin{proof} Applying \eqref{green-project} to the projection $\pi_1$ gives
$G_1(x_1,y_1) = \sum_{w_2 \in H(y_2)} G(x,y_1w_2)\,.$

Let $w_2 \in H(y_2)$, where $H(y_2)$ is the horocycle of $y_2$ in $\T_r$. 
We write $v_2 = v(w_2)$ for the unique element in $H(x_2)$ that satifies 
$v_2 \lle w_2$. 
By Lemma \ref{stop}, the random walk has to pass through some point of
the form in $\{u_1v_2 :u_1 \in H(x_1)\}$ on the way from $x$ to $y_1w_2$, 
that is, 
$$
\begin{aligned}
G(x,y_1w_2) &= \Ex_{x}\Bigl(\1_{[\tm_2(v_2) < \infty]}\, 
G\bigl(Z_{\tm_2(v_2)},y_1w_2\bigr) \Bigr)\\
&= \Ex_{x}\Bigl(\1_{[\tm_2(v_2) < \tm_1(c_1)]}\, 
G\bigl(Z_{\tm_2(v_2)},y_1w_2\bigr) \Bigr)
+ \Ex_{x}\Bigl(\1_{[\tm_1(c_1) < \tm_2(v_2) < \infty]}\, 
G\bigl(Z_{\tm_2(v_2)},y_1w_2\bigr) \Bigr)\,.
\end{aligned}
$$
Now, if starting at $x$, we have 
$\tm_2(v_2) < \tm_1(c_1)$, then $Z_{\tm_2(v_2)} = u_1v_2$ for some random
$u_1 \in H(x_1)$ that must satisfy $\up(u_1,y_1)=\up_1$ and 
$\dn(u_1,y_1)=\dn_1$, since $c_1$ cannot lie on $\geo{x_1\,u_1}$. 
But we also have $\up(v_2,w_2) = \up_2=0$ and
$\dn(v_2,w_2) = \dn_2$. That is, the points $u_1v_2$ and $y_1w_2$ habe the
same relative position as the points $x$ and $y$, and therefore
$G(u_1v_2,y_1w_2) = G(x,y)$ by Lemma \ref{coordinates}. We get
$$
\Ex_{x}\Bigl(\1_{[\tm_2(v_2) < \tm_1(c_1)]}\, 
G\bigl(Z_{\tm_2(v_2)},y_1w_2\bigr) \Bigr)
= \Prb_x[\tm_2(v_2) < \tm_1(c_1)] \,G(x,y)\,.
$$
Now, given $v_2 \in H(x_2)$, there are precisely $r^{\dn_2}$ elements 
$w_2 \in H(y_2)$ with $v(w_2)=v_2$. Combining all these observations, 
$$
\begin{aligned}
G_1(x_1,y_1) &= 
\Biggl(\sum_{v_2 \in H(x_2)} \Prb_x[\tm_2(v_2) < \tm_1(c_1)]\Biggr) 
r^{\dn_2}\,G(x,y) + R(\dn_1,\dn_2)\,,\quad\mbox{where}\\
R(\dn_1,\dn_2) &= \sum_{w_2 \in H(y_2)} 
\Ex_{x}\Bigl(\1_{[\tm_1(c_1) < \tm_2(v(w_2)) < \infty]}\, 
G\bigl(Z_{\tm_2(v(w_2))},y_1w_2\bigr) \Bigr)\,.
\end{aligned}
$$
Let us first consider the error term $R(\dn_1,\dn_2)$. Note that 
$G(\cdot,\cdot) \le G(o,o) < \infty\,$, since our random walk is transient.
(already the projections onto $\T_q$ and $\T_r$ are transient!)
Since $\dn_2 \to \infty$, also $\tm_1(c_1) \to \infty$ almost surely.
It follows that
$$
\begin{aligned} 
r_{\dn_1,\dn_2}(w_2)
&:= \Ex_{x}\Bigl(\1_{[\tm_1(c_1) < \tm_2(v(w_2)) < \infty]}\, 
G\bigl(Z_{\tm_2(v(w_2))},y_1w_2\bigr) \Bigr) \\
&\le \Prb_x[\tm_1(c_1) < \tm_2(v(w_2))< \infty]\,G(o,o) \to 0 \quad\mbox{when}\;
\dn_2 \to \infty\,.
\end{aligned}
$$ 
On the other hand,
$$
r_{\dn_1,\dn_2}(w_2) \le G(x,y_1w_2) \AND 
\sum_{w_2 \in H(y_2)} G(x,y_1w_2) = G_1(x_1,y_1) \le G_1(o_1,o_1)\,.
$$
Thus, dominated convergence (in the summation) implies that
$R(\dn_1,\dn_2) \to 0$ as $\dn_2 \to \infty$.

\smallskip

It remains to show that 
$C(\spn) = \sum_{v_2 \in H(x_2)} \Prb_x[\tm_2(v_2) < \tm_1(c_1)]$
has the proposed asymptotic behaviour, when $\dn_2$ (and $\spn$) $\to \infty$.

We may suppose without loss of generailty that $\hor(x_1) = \hor(x_2) = 0$,
so that $\hor(c_1) = -\spn$. Then Lemma \ref{stop} implies 
$\Prb_x[\tm_2(v_2) < \tm_1(c_1)] = \Prb_x[\tm_2(v_2) < \wt \tm(-\spn)]$.
Now let the superscript $^{(-\spn)}$ refer to the random walk $\wt P$ on $\Z$
stopped at $-\spn$, i.e., we consider the restriction of $\wt P$
to $\{ k \in \Z : k > -\spn \}$. We use the same superscript for the
random walk $P_2$ on $\T_r$ stopped at the horocycle $H_{\spn}$ in $\T_r$,
i.e., we consider the restriction of $P_2$ to 
$\{ z_2 \in \T_r : \hor(z_2) < \spn \}$.
Using the notation of \eqref{G1xy}, we have 
$$
\Prb_x[\tm_2(v_2) < \wt \tm(-\spn)] = F_2^{(-\spn)}(x_2,v_2) = 
\frac{G_2^{(-\spn)}(x_2,v_2)}{G_2^{(-\spn)}(v_2,v_2)} =
\frac{G_2^{(-\spn)}(x_2,v_2)}{G_2^{(-\spn)}(o_2,o_2)}\,.
$$
Since $\wt P^{(-\spn)}$ is the projection of $P_2{(-\spn)}$ under the
mapping $z_2 \mapsto -\hor(z_2)$, we find
$$
C(\spn) = \sum_{v_2 \in H(x_2)} 
\frac{G_2^{(-\spn)}(x_2,v_2)}{G_2^{(-\spn)}(o_2,o_2)} = 
\frac{\wt G^{(-\spn)}(0,0)}{G_2^{(-\spn)}(o_2,o_2)}\,.
$$
If $\spn \to \infty$ then $G_2^{(-\spn)}(o_2,o_2) \to G_2(o_2,o_2) < \infty$,
for each value of $\al$. 
If $\al < 1/2$ then $G^{(-\spn)}(0,0) \to \wt G(0,0) < \infty$. 
If $\al=1/2$ then routine calculations regarding SRW on $\Z$ yield
$G^{(-\spn)}(0,0) = 2\spn$. Thus, $C(\spn)$ has the proposed
asymptotic behaviour.
\end{proof}

The last proposition is valid for arbitrary $\al$. However, it becomes
meaningful only when $\al \le 1/2$. Indeed, when $\al > 1/2$, then the left hand
side in the decomposition \eqref{vertical} tends to $0$ by \eqref{G1xy}.
In this case, \eqref{vertical} contains no information about the
asymptotic behaviour of $G(x,y)$. On the other hand, when $\al \le 1/2$
and $\dn_1 = 0$ (situation (II) of Figure 4) then $G_1(x_1,y_1) = G_1(o_1,o_1)$
is constant, see \eqref{G1xy}. When we consider the ``dual'' situation of 
Figure 5, as illustrated in Figure 6, this discussion shows that it is
not useful to rewrite Proposition \ref{IIbis} by just exchanging \emph{both}
the roles of the two trees \emph{and} $\al$ with $1-\al$. 

We shall use the superscript $^*$ for the respective random walks on
$\DL$, $\T_q$, $\T_r$, and $\Z$ that are obtained by exchanging
$\al \leftrightarrow 1-\al$, \emph{without} exchanging roles of the two trees.
Thus, $P_{\al}^* = P_{1-\al}\,$, $P_{\al,q}^*=P_{1-\al,q}\,$, 
$P_{1-\al,r}^*=P_{\al,r}\,$, and $\tilde P^*$ on $\Z$ moves from $k$ to $k+1$
with probability $1-\al$ and to $k-1$ with probability $\al$.

\begin{lem}\label{exchange}
$\;G^*(x,y)= \be^{\hor(y_1)-\hor(x_1)}\,G(x,y) \quad \forall\ x,y \in \DL
\quad (\be = \frac{1-\al}{\al}).$
\end{lem}

\begin{proof} The function $g(x) = \be^{\hor(x_1)}$ satisfies $Pg=g$, and
$p^*(x,y) = p(x,y)\,g(y)/g(x).$
\end{proof} 
$$
\beginpicture
\setcoordinatesystem units <2.2mm,2.56mm> 

\setplotarea x from -24 to 4, y from -7.5 to 3

\put{$y_2$} [lt] at -1.85 -2.2
\put{$c_2$} [rt] at -4.15 -4.2
\put{$\omega_2$} [lt] at -0.7 -7.7 
\put{$x_2$} [lb] at -12 4.2

\arrow <5pt> [.2,.67] from -12 4 to -0.5 -7.5 
\plot -4 -4  -2 -2 /

\put{$x_1$} [lb] at  -15.95 4.2
\put{$\omega_1$} [rb] at -19.7 7.7
\put{$y_1$} [rt] at  -9.8 -2.2

\arrow <5pt> [.2,.67] from -10 -2 to -19.5 7.5

\multiput {$\scs\bullet$} at -2 -2  -4 -4  -12 4  -10 -2  -16 4  /

\setdots<4pt>

\plot -20.5 -2  -2 -2 /
\plot -25 -4  -4 -4 /
\plot -12 4  -25 4 /

\put{$\spn= \up_2\!\!\left\{\rule[11mm]{0mm}{0mm}\right.$}[rb] at -25 -3.95
\put{$\dn_2\!\left\{\rule[3.5mm]{0mm}{0mm}\right.$}[rt] at -20.5 -2
\put{$\dn_1\!\!\left\{\rule[-7.34mm]{0mm}{0mm}\right.$}[rb] at -20.5 -2.25

\endpicture
$$
\begin{center}
\emph{Figure 6}
\end{center}

\begin{cor}\label{IIIbis} If, as in Figure 6, $x=x_1x_2$ and $y=y_1y_2$ 
with $\spn = \spn(x,y)$ are such that $\dn_1 = \up_2 - \dn_2 \to \infty$, 
$\up_1 = 0$ and $\dn_2$ is arbitrary, then
\begin{equation}\label{vertical*}
G_2^*(x_2,y_2) = C^*(\spn)  \, (q\,\be)^{\dn_1}\,G(x,y) 
+ R^*(\dn_1,\dn_2)\,,
\end{equation}
where
$$
C^*(\spn) \to \frac{\wt G(0,0)}{G_1(o_1,o_1)}\quad
\mbox{if}\;\;\al \ne \frac12\,,\qquad\quad
C^*(\spn)  \sim \frac{2\spn}{G_1(o_1,o_1)}\quad
\mbox{if}\;\;\al = \frac12\,,
$$
and $0 < R^*(\dn_1,\dn_2) < G_2^*(x_2,y_2)$ with 
$$
\lim_{\dn_1 \to \infty} R^*(\dn_1,\dn_2) = 0\,.
$$
\end{cor}

This is immediate by applying Proposition \ref{IIbis} to $P^*$ with
$r \lra q$. Also observe that $\wt G^*(0,0)= \wt G(0,0)$ and  
$G_i^*(o_i,o_i) = G_i(o_i,o_i)$ for $i=1,2$.
In the specific case $\dn_1=0$ (resp. $\dn_2=0$), Proposition \ref{IIbis}
(resp. Corollary \ref{IIIbis}) yields the asymptotic behaviour of
$G(x,y)$ in situation (II) (resp. (III)) of Figure 4.

\begin{cor}\label{II-III} {\rm (a)} Referring to situation (II) of Figure 4,
if $\spn = \up_1 = \dn_2 \to \infty$ and $\dn_1=\up_2=0$ then
$$
G(x,y) \sim \begin{cases}
\dfrac{G_1(o_1,o_1)\,G_2(o_2,o_2)}{\wt G(0,0)\, r^{\spn}} 
&\mbox{if}\;\;\al > 1/2\,,\AND\\[15pt]
\dfrac{G_1(o_1,o_1)\,G_2(o_2,o_2)}{2\,\spn\, r^{\spn}}
&\mbox{if}\;\;\al = 1/2\,.
\end{cases}
$$
{\rm (b)} Referring to situation (III) of Figure 4,
if $\spn = \dn_1 = \up_2 \to \infty$ and $\up_1=\dn_2=0$ then
$$
G(x,y) \sim \begin{cases}
\dfrac{G_1(o_1,o_1)\,G_2(o_2,o_2)}{\wt G(0,0)\, (q\,\be)^{\spn}} 
&\mbox{if}\;\;\al > 1/2\,,\AND\\[15pt]
\dfrac{G_1(o_1,o_1)\,G_2(o_2,o_2)}{2\,\spn\, q^{\spn}}
&\mbox{if}\;\;\al = 1/2\,.
\end{cases}
$$
\end{cor}

Proposition \ref{IIbis}, resp. Corollary \ref{IIIbis}, also leads to 
an asymptotic estimate of $G(x,y)$ when $\dn_1$, resp. $\dn_2$, remains
bounded. Otherwise, the left hand side of the decomposition \eqref{vertical},
resp. \eqref{vertical*}, tends to $0$. Nevertheless, those decompositions
will now be useful ``on the average'' for situations (I) and (IV) of
Figure 4. 

\begin{pro}\label{IV} Referring to situation (IV) of Figure 4, if
$x=x_1x_2$ and $y=y_1y_2$ with $\spn = \spn(x,y)$ are such that 
$\up_1=\dn_1=\spn \to \infty$ and $\up_2=\dn_2=0$ ($y_2=x_2$) then
$$
G(x,y) \sim \begin{cases}
B^*_{\be}\,\dfrac{G_1(o_1,o_1)\,G_2(o_2,o_2)}{\wt G(0,0)\,(q\,\be)^{\spn}}\,,
&\mbox{where}\quad B^*_{\be}= \dfrac{(\be-1)(r\,\be+1)}{r\,\be^2-1}\,,
\quad\mbox{if}\;\ \al < 1/2\,,\\[12pt]
B^*_1\,\dfrac{G_1(o_1,o_1)\,G_2(o_2,o_2)}{2\,\spn^2\,q^{\spn}}\,,
&\mbox{where}\quad B^*_1= \dfrac{r+1}{r-1}\,,
\quad\mbox{if}\;\ \al = 1/2\,.
\end{cases}
$$
\end{pro}

\begin{proof} Again, we may assume that $\hor(x_1)=\hor(x_2)=0$.
Since $1-\al \ge 1/2$, 
$$
\tm := \tm_1(c_1) = \wt \tm(-\spn) < \infty \quad \Prb_x\mbox{-almost surely.}
$$
This and Lemma \ref{stop} yield
$$
G(x,y) = 
\Ex_x\bigl( \1_{[\tm < \infty]}\,G(Z_{\tm},y)\bigr)
 = \Ex_x\bigl(G(Z_{\tm},y)\bigr)\,.
$$
We have $\pi_1(Z_{\tm}) = c_1$. Set $D_{\tm} = \dn(Z_{\tm}^2,x_2)$,
a non-negative, integer-valued random variable, see Figure 7. 
$$
\beginpicture
\setcoordinatesystem units <2.5mm,3.1mm> 

\setplotarea x from -4 to 24, y from -3.8 to 8.5
\arrow <5pt> [.2,.67] from 4 4 to 0 8 
\plot 8 0  4 4  0 0 /

\plot 12 0   16 -4   24 4 /         
\arrow <5pt> [.2,.67] from 16 -4 to 20 -8

\put{$x_1$} [lb] at  8.15 0.2
\put{$y_1$} [rb] at -.15 0.2
\put{$c_1$} [lb] at 4.15 4.2
\put{$\omega_1$} [rb] at 0.2 8.2

\put{$x_2=y_2$} [lb] at  12.15 0.2
\put{$\omega_2$} [lt] at 20.2 -8.2
\put{$Z_{\tm}^2$} [lb] at  24.15 4.2

\multiput {$\scs\bullet$} at 
                            0 0  4 4  8 0  12 0  16 -4  24 4   8 4 /

\setdashes <3pt>

\plot  12 0   8 4 /


\setdots<4pt>
\plot -9.5 0  12 0 /
\plot -9.5 4  24 4 /
\plot -9.5 -4  16 -4 /

\put{$\spn\!\left\{\rule[-5mm]{0mm}{0mm}\right.$}[rb] at -7.5 -.1 
\put{$D_{\tm}\!\!\left\{\rule[-5mm]{0mm}{0mm}\right.$}[rt] at -7.5 0

\endpicture
$$


\begin{center}
\emph{Figure 7}
\end{center}
The relative position of $c_1 Z_{\tm}^2$ (in the place of $x$) and $y=y_1x_2$ 
is precisely the one of Figure 6, replacing $\spn=\spn(x,y)$ with 
$\spn+D_{\tm}$ and $\dn_2$ with $D_{\tm}$. We can apply Corollary \ref{IIIbis}.
Note that by \eqref{G1xy}, applied to $P_2^* = P_{\al,r}$,
$$
G_2^*(Z_{\tm}^2,x_2) = (r\,\be)^{-D_{\tm}} G_2^*(x_2,x_2) = 
(r\,\be)^{-D_{\tm}} G_2(o_2,o_2)\,.
$$
We get
\begin{equation}\label{GZy}
\begin{aligned}
G(Z_{\tm},y) &= \frac{1}{C^*(\spn+D_{\tm})\,(q\,\be)^{\spn}}
\Bigl( G_2^*(Z_{\tm}^2,x_2) - R^*(\spn, D_{\tm}) \Bigr)\\
&= \frac{1}{C^*(\spn+D_{\tm})\,(q\,\be)^{\spn}}
\left( \frac{1}{(r\,\be)^{D_{\tm}}} G_2(o_2,o_2)\ - R^*(\spn, D_{\tm})
\right)\,.
\end{aligned}
\end{equation}\\
\underline{Case 1.} $\al < 1/2$. Then $C^*(\spn+D_{\tm}) \to \tilde
G(0,0)/G_1(o_1,o_1)\,$, a finite limit. Since 
$R^*(\spn,D_{\tm}) < G_2^*(Z_{\tm}^2,x_2) \le G_2^*(o_2,o_2)$ and
$R^*(\spn,D_{\tm}) \to 0$ as $\spn \to \infty$, dominated convergence yields
\begin{equation}\label{rest}
\Ex_x[R^*(\spn, D_{\tm})] \to 0\,.
\end{equation}
Also, $Z_n^2$ converges almost
surely to a $\bd^*\T_r$-valued random variable $Z_{\infty}^2$; see
{\sc Cartwright, Kaimanovich and Woess} \cite{CaKaWo}, where this is proved
under much more general assumptions. Since $\tm \to \infty$ when 
$\spn \to \infty$, we get $Z_{\tm}^2 \to Z_{\infty}^2$ and consequently
$$
D_{\tm} \to D_{\infty} = d(o_2, o_2 \cf Z_2^{\infty})\,.
$$
(Cf. \S \ref{geometry} for notation.) 
Therefore, \eqref{GZy} and \eqref{rest} yield
$$
\Ex_x\bigl(G(Z_{\tm},y)\bigr) \sim 
\frac{G_1(o_1,o_1)G_2(o_2,o_2)}{\wt G(0,0)\,(q\,\be)^{\spn}}\,B^*_{\be}\quad
\mbox{as}\; \spn \to \infty\,,
$$
where $B^*_{\be} = \Ex_x\bigl[(r\,\be)^{-D_{\infty}}\bigr]\,$. This number can be
computed explicitly: let $w_2^{(m)}$ denote the vertex on $\geo{x_2\,\om_2}$
at distance $m$ from $x_2$. If $m \ge 1$, then $D_{\infty} \ge m$
precisely when $Z_{\infty}^2 \in \wh \T_r(x_2,w_2^{(m)})$. Applying a
frequently used formula for the limit distribution on the boundary of
arbitrary transient nearest neighbour random walks on trees (see e.g.
\cite{Car}), we get that
$$
\Prb_x[D_{\infty} \ge m] = 
\frac{F_2(x_2,w_2^{(m)})\bigl(1-F_2(w_2^{(m)},w_2^{(m-1)})\bigr)}
     {1 - F_2(w_2^{(m-1)},w_2^{(m)})F_2(w_2^{(m)},w_2^{(m-1)})}
= \be^{-m}\frac{\be\,r - \be}{\be\,r - 1}\,.
$$
We have used the $P_2$-version of \eqref{G1xy} in the last computation.
It is now straghtforward that $B^*_{\be}$ has the proposed value.\\[4pt]
\underline{Case 2.} $\al = 1/2$. 
Here, we need to compute explicitly the distribution of $D_{\tm}$, 
which depends on $\spn$. 
Consider the random variable $M= M_{\spn} = \max \{ \wt Z_n : n < \tm \}$.
If $\nn = \nn_{\spn} = \max \{ n < \tm : \wt Z_n = M_{\spn}\}$, then
$Z_{\nn}^2$ must be the point $w_2 = w_2^{(M)}$ on $\geo{o_2\,\om_2}$. 

Conditioned on the value of $M_{\spn}$,
the random element $Z_{\tm}^2$ is equidistributed on the set
$\{ v_2 \in \T_r : \hor(v_2) = \spn\,,\;w_2 \lle v_2 \}$, 
which has $r^{\spn+M}$ elements. Among the latter, the number of elements
with $\dn(v_2,x_2) = d \in \{0,\dots, M\}$ is $r^{\spn}$, if $d=0$,
and $(r-1)r^{\spn+d-1}$, if $d \ge 1$. If $D_{\tm} = d$ then 
$M_{\spn} \ge d$. Thus, if $0 \le d \le m$ then
$$
\Prb_x[D_{\tm} = d \mid M=m] = \epsilon_d \,r^{d-m} \,,\quad \mbox{where}
\quad \epsilon_d = 
 \begin{cases} 1\,&\mbox{if}\;d=0\,,\\
                (r-1)/r\,,&\mbox{if}\;d \ge 1\,.
\end{cases}
$$
Also, $\Prb_x[M \ge m]$ is the probablity that the random walk $\wt Z_n$		
on $\Z$ reaches $m$ before $-\spn$. This is 
$\varphi_2(\spn,m) = \spn/(\spn+m)$, as given in
\eqref{varphi}, and $\Prb_x[M = m] = \varphi_2(\spn,m)-\varphi_2(\spn,m+1)$.
We find
\begin{equation}\label{D-dist}
\Prb[D_{\tm} = d] = 
\epsilon_d  \sum_{m=d}^{\infty} \frac{r^d}{r^{m}} \,
\frac{\spn}{(\spn+m)(\spn+m+1)}\,.
\end{equation}
We know from Corollary \ref{IIIbis} that $C^*(\spn) \sim 2\spn/G_1(o_1,o_1)$.
Therefore \eqref{GZy} implies
$$
\begin{gathered}
2\,\spn^2\, q^{\spn}\,G(x,y) \sim 
G_1(o_1,o_1)G_2(o_2,o_2) \,
\Ex_x\left( \frac{\spn^2}{(\spn+D_{\tm})\,r^{D_{\tm}}}\right) - \Rest(\spn)\,,\\
\mbox{where}\qquad
\Rest(\spn) =\Ex_x
\left(\frac{\spn^2}{\spn+D_{\tm}}\,R^*(\spn,D_{\tm})\right).
\end{gathered}
$$
Using \eqref{D-dist}, we can write
$$
\begin{gathered}
\Ex_x\left( \frac{\spn^2}{(\spn+D_{\tm})\,r^{D_{\tm}}}\right) =
\sum_{d=0}^{\infty} f_{\spn}(d)\qquad\mbox{with}\\
f_{\spn}(d) = \Prb_x[D_{\tm}=d]\,\frac{\spn^2}{(\spn+d)\,r^{d}} 
= \epsilon_d \,\frac{\spn}{\spn+d} 
\sum_{m=d}^{\infty} r^{-m} \frac{\spn^2}{(\spn+m)(\spn+m+1)}\,.
\end{gathered}
$$
Now $f_{\spn}(d)$ is increasing in $\spn$, and 
$f_{\spn}(d) \to f(d)$ with $f(0) = r/r-1$ and $f(d) = r^{-d}$ 
for $d \ge 1$. Monotone convergence implies that
$$
\sum_{d=0}^{\infty} f_{\spn}(d) \to \sum_{d=0}^{\infty} f(d) = 
\frac{r+1}{r-1}=B^*_1\,.
$$
To conclude our asymptotic estimate, we have to show that the rest
tends to zero. We expand
$$
\Rest(\spn) 
= \sum_{d=0}^{\infty} \Prb_x[D_{\tm}=d]\,\frac{\spn^2}{\spn+d} \,R^*(\spn,d) \,
= \sum_{d=0}^{\infty} f_{\spn}(d)\, r^d\,R^*(\spn,d)\,.
$$
We have 
$
R^*(\spn,D_{\tm}) < G^*_2(Z_{\tm}^2,x_2) = r^{-D_{\tm}}\,G_2(o_2,o_2)\,.
$
and $f_{\spn}(d)\, r^d\,R^*(\spn,d) < f(d)\, G_2(o_2,o_2)$.
Also,  $\sum_d f(d) < \infty$. On the other hand,  
$\lim_{\spn \to \infty} f_{\spn}(d)\, r^d\,R^*(\spn,d) = 0$ pointwise in $d$. 
Dominated convergence implies $\Rest(\spn) \to 0$.
\end{proof}

We remark that in the proof we might have treated Case 1 in the same way 
as Case 2, by first determining the distribution of $D_{\tm}$ and then letting 
$\spn \to \infty$ (whence $\tm \to \infty$). However, it is more likely that
the method used above will lend itself to an extension to finite range
(instead of nearest neighbour) random walks on $\DL(q,r)$ where $p(x,y)$
depends only on $\up_1, \dn_1, \up_2, \dn_2$.

Again, from the last proposition we can also deduce the asymptotics of
$G(x,y)$ in the dual situation (I) of Figure 4 by considering $P^*$. 
In this case we have
$G^*(x,y) = G(x,y)$, since $\hor(x_1) - \hor(y_1) = 0$. Thus, we only have to
exchange $r \leftrightarrow q$.

\begin{cor}\label{I} Referring to situation (I) of Figure 4, if
$x=x_1x_2$ and $y=y_1y_2$ with $\spn = \spn(x,y)$ are such that 
$\up_1=\dn_1=0$ ($y_1=x_1$) and $\up_2=\dn_2=\spn \to \infty$ then
$$
G(x,y) \sim \begin{cases}
B_{\be}\,\dfrac{G_1(o_1,o_1)\,G_2(o_2,o_2)}{\wt G(0,0)\,(r\,\be)^{\spn}}\,,
&\mbox{where}\quad B_{\be}= \dfrac{(\be-1)(q\,\be+1)}{q\,\be^2-1}\,,
\quad\mbox{if}\;\ \al < 1/2\,,\\[12pt]
B_1\,\dfrac{G_1(o_1,o_1)\,G_2(o_2,o_2)}{2\,\spn^2\,r^{\spn}}\,,
&\mbox{where}\quad B_1= \dfrac{q+1}{q-1}\,,
\quad\mbox{if}\;\ \al = 1/2\,.
\end{cases}
$$
\end{cor}
 
Combining Corollary \ref{II-III}, Proposition \ref{IV} and Corollary \ref{I}
with Proposition \ref{decomp}, we obtain Theorem \ref{estimates}.

\section{The Martin compactification}\label{martin}

We are now ready to determine the full Martin compactifiction of
$P=P_{\al}$ on $\DL(q,r)$. Recall that the Martin compactification 
of the projected random walk $P_{\al,q}$ is $\wh \T_q$, the end
compactification of the tree. (The analogous result holds of course
for the second projection $P_{1-\al,r}$ on $\T_r$.) 
The end compactification of $\T_q$ was described in \S \ref{geometry};
in particular, it is a compact metric space with the ultrametric
\begin{equation}\label{theta}
\theta(z_1,w_1) = \begin{cases} 0\,,&\mbox{if}\quad z_1=w_1\,,\\
q^{d(c_1,o_1)}\,,&\mbox{where} \; c_1 = z_1 \wedge w_1 \;
(\mbox{confluent w.r.t.}\; o_1)\,,\quad\mbox{if}\quad z_1 \ne w_1
\end{cases}
\end{equation}
for $z_1, w_1 \in \wh \T_q$. In particular (recall), 
$z_1^{(n)} \to \xi_1 \in \bd\T_q$ if and only if 
$d(z_1^{(n)} \wedge \xi_1,o_1) \to \infty$.

Since $\DL(q,r) \subset \T_q \times \T_r\,$,
this provides us with a natural \emph{geometric compactification}
$\wh{\DL}(q,r)$, namely, the closure of $\DL(q,r)$ in 
$\wh\T_q \times \wh\T_r$. The ideal boundary of $\DL$ in this compactification
consists of $5$ disjoint pieces:
\begin{equation}\label{bdry-I}
\bigl(\bd^*\T_q \times \{ \om_2 \}\bigr) \cup
\bigl(\{ \om_1 \} \times \bd^*\T_r\bigr) \cup
\bigl\{ \om_1 \om_2 \bigr\} \cup
\bigl(\T_q \times \{ \om_2 \}\bigr) \cup
\bigl(\{ \om_1 \} \times\T_r\bigr)\,,
\end{equation}
compare with \cite{Ber}. For a sequence $y^{(n)} = y_1^{(n)}y_2^{(n)} \in \DL$,
we have
\begin{equation}\label{convergence-I}
\begin{alignedat}{3}
y^{(n)} &\to \xi_1\om_2\,,\quad \xi_1 \in \bd^*\T_q\,,
&\quad&\mbox{if }\, y_1^{(n)} \to \xi_1 &\AND y_2^{(n)} &\to \om_2\,;\\
y^{(n)} &\to \om_1\xi_2 \,,\quad \xi_2\in \bd^*\T_r\,,
&\quad&\mbox{if}\quad y_1^{(n)} \to \om_1 &\AND y_2^{(n)} &\to \xi_2\,;\\
y^{(n)} &\to \om_1\om_2 \,,
&\quad&\mbox{if}\quad y_1^{(n)} \to \om_1 &\AND y_2^{(n)} &\to \om_2\,;\\
y^{(n)} &\to y_1\om_2\,,\quad y_1\in \T_q \,,
&\quad&\mbox{if}\quad y_1^{(n)}=y_1\; \forall\ n \ge n_0 &\AND y_2^{(n)} &\to \om_2\,;\\
y^{(n)} &\to \om_1y_2\,,\quad y_2 \in \T_r\,,
&\quad&\mbox{if}\quad y_1^{(n)}\to \om_1 &\AND y_2^{(n)}&=y_2\; \forall\ n \ge n_0\,.
\end{alignedat}
\end{equation}
Every sequence in $\DL$ that tends to infinity has a subsequence of
one of these 5 types. 

Recall from \S \ref{trees} the Martin kernels associated with $P_{\al,q}$
and $P_{1-\al,r}$ and the spherical functions \eqref{spherical}.

\begin{thm}\label{centred}
If $\al =1/2$ then the Martin compactification of $\DL(q,r)$
with respect to $P = P_{1/2}$ is the geometric compactification
$\wh{\DL}(q,r)$. The extension of the Martin kernel on the boundary
described in \eqref{bdry-I} and \eqref{convergence-I} is given by
\begin{alignat}{2}
K(x_1x_2, \xi_1\om_2) &= K_1(x_1,\xi_1)\,,&\quad &\xi_1 \in
\bd^*\T_q\,,\tag{\rm i}\\
K(x_1x_2, \om_1\xi_2) &= K_2(x_1,\xi_1)\,,&\quad &\xi_2\in \bd^*\T_r\,,
\tag{\rm ii}\\
K(x_1x_2, \om_1\om_2)  &=1\,,&\quad&\tag{\rm iii}\\
K(x_1x_2, y_1\om_2) &= \frac{\phi_{\frac12,q}(x_1,y_1)}
{\phi_{\frac12,q}(o_1,y_1)}\,,
&\quad& y_1\in \T_q \,, \AND\tag{\rm iv}\\
K(x_1x_2, \om_1y_2) &= \frac{\phi_{\frac12,r}(x_2,y_2)}
{\phi_{\frac12,r}(o_2,y_2)}\,,
&\quad& y_2 \in \T_r\,.\tag{\rm v}
\end{alignat}
Each of the kernels in\/ {\rm (i), (ii), (iii)} constitutes a minimal harmonic
function, while the ones of\/ {\rm (iv)} and\/ {\rm (v)} are non-minimal
harmonic. 
\end{thm}

\begin{proof} We use part (b) of Theorem \ref{estimates}.

\smallskip

(i) Suppose first that $y = y_1y_2 \to \xi_1\om_2$. Then, with the usual
notation $\up_i=\up(x_i,y_i)$ and $\dn_i=\dn(x_i,y_i)$, we see
that $\up_1 = d(x_1, x_1 \wedge \xi_1)$ is constant when $y_1$ is sufficiently
close to $\xi_1$. On the other hand, $\dn_1, \up_2-\dn_2 \to \infty$. 
Recall that
$\up_1 + \up_2 = \dn_1 + \dn_2 = \spn \to \infty$. The dominant term in
$$
\frac{q+1}{q-1}\,\up_2\,\dn_2 +\spn\,\up_2\,\dn_1 + 
\spn\,\up_1\,\dn_2 + \frac{r+1}{r-1}\,\up_1\,\dn_1 
$$
is $\spn\,\up_2\,\dn_1$, since 
$\frac{\up_2\,\dn_2}{\spn\,\up_2\,\dn_1} \le \frac{1}{\dn_1}\,$,
$\;\frac{\spn\,\up_1\,\dn_2}{\spn\,\up_2\,\dn_1} \le \frac{\up_1}{\dn_1}\,$
and $\frac{\up_1\,\dn_1}{\spn\,\up_2\,\dn_1} \le \frac{\up_1}{\spn}$
all tend to zero. As $\up_2 \sim \spn$, we find
$$
G(x,y) \sim A_1 \frac{\dn_1}{\spn^2\,q^{\dn_1}\,r^{\dn_2}} =
A_1 \frac{\dn(x_1,y_1)}{\spn(x,y)^2}\,F_1(x_1,y_1)\,F_2(x_2,y_2)\,.
$$
For each $x \in \DL$ we have that $\spn(x,y)-\spn(o,y)$ and 
$\dn(x_1,y_1) - \dn(o_1,y_1)$ are constant when
$y_1$ is close to $\xi_1$ in the end metric $\theta$. 
Therefore
$$
K(x,y) \sim \frac{\spn(o,y)^2\, \dn(x_1,y_1)}{\spn(x,y)^2\, \dn(o_1,y_1)}
\,K_1(x_1,y_1)\,K_2(x_2,y_2) \to K_1(x_1,\xi_1)\,K_2(x_2,\om_2) 
= K_1(x_1,\xi_1)\,,
$$
since $K_2(\cdot,\om_2) \equiv 1$.

\smallskip

(ii) follows immediately from (i), exchanging $r \leftrightarrow q$.

\smallskip

(iii) If $y=y_1y_2 \to \om_1\om_2$ then $\up_i=\up(x_i,y_i) \to \infty$ 
for $i=1,2$. For given $x$, when $\up_i > \up(x_i,o_i)$ for $i=1,2$, then
$\dn_i = \dn(x_i,y_i)$ coincides with $\dn(o_i,y_i)$ and
$\spn = \spn(x,y)$ coincides with $\spn(o,y)$, and we also have 
$\up_1 - \up(o_1,y_1) = \up(o_2,y_2) - \up_2 = k$, where $k = \hor(x_1)$.
Therefore
$$
K(x,y) \sim \frac{\frac{q+1}{q-1}\,\up_2\,\dn_2 +\spn\,\up_2\,\dn_1 + 
\spn\,\up_1\,\dn_2 + \frac{r+1}{r-1}\,\up_1\,\dn_1}
{\frac{q+1}{q-1}\,(\up_2+k)\,\dn_2 +\spn\,(\up_2+k)\,\dn_1 + 
\spn\,(\up_1-k)\,\dn_2 + \frac{r+1}{r-1}\,(\up_1-k)\,\dn_1} \to 1\,.
$$

(iv) If $y_1$ remains fixed and $y_2 \to \om_2$, then
$\dn_1 = \dn(x_1,y_1)$ and $\up_1 = \up(x_1,y_1)$ are constant.
Since $\up_2 = \spn - \up_1$ and $\dn_2 = \spn-\dn_1$, where
$\spn=\spn(x_1,y_1)$, we get
$$
\begin{aligned}
G(x,y) &\sim \frac{A_1}{\spn^2\,q^{\dn_1}\,r^{\dn_2}}
\left(\frac{q+1}{q-1}\!
\left(1\!-\!\frac{\up_1}{\spn}\right)\!\!\Bigl(1\!-\!\frac{\dn_1}{\spn}\Bigr)
+\left(1\!-\!\frac{\up_1}{\spn}\right)\dn_1 + 
\up_1\Bigl(1\!-\!\frac{\dn_1}{\spn}\Bigr)+ 
\frac{r+1}{r-1}\,\frac{\up_1\,\dn_1}{\spn^2} \right) \\
&\sim \frac{A_1}{\spn^2\,q^{\dn_1}\,r^{\dn_2}}
\left(\frac{q+1}{q-1} + \dn_1 + \up_1\right) 
=  \,\frac{A_1}{\spn^2\,r^{\dn_2}}\,\phi_{\frac12,q}(x_1,y_1)\,.
\end{aligned}
$$
As above in (iii), $\dn(o_2,y_2) = \dn_2$ ($=\dn(x_2,y_2)$) when 
$\up(x_2,y_2) > \up(x_2,o_2)$, and then also $\spn(o_1,y_1) - \spn = k$, a
constant. Therefore, 
$$
G(o,y) \sim \frac{A_1}{(\spn+k)^2\,r^{\dn_2}}\,\phi_{\frac12,q}(o_1,y_1)\,.
$$
Thus, we obtain the proposed limit of $K(x,y)$ as $y \to y_1\om_2)$.

\smallskip

(v) follows from (iv), exchanging $r \leftrightarrow q$.

\smallskip

Finally, the -- simple -- proof of minimality of the functions in 
(i), (ii) and (iii) can be found in \cite{Woe}. Non-minimality of the
spherical functions in (iv) and (v) is straightforward, since they 
are also non-minimal for the projected random walks on the respective trees. 
\end{proof}

Next, we explain what happens in the case $\al \ne 1/2$.
If $(y_1^{(n)})$ is a sequence in $\T_q$ with $\up(o_1,y_1^{(n)}) \to \infty$
then $y_1^{(n)} \to \om_1\,$, independently of the values of $\hor(y_1^{(n)})$.
The horocyclic drawing of $\T_q$ as in Figure 1 suggests that one
may use a finer distinction by introducing boundary points 
$\om_1^k$, $k \in \overline\Z = \Z \cup \{\pm \infty\}$, at 
infinity, one for each horocycle, one at the ``level'' $-\infty$, and one at 
the level $+\infty$ (thinking of $\bd^*\T_q$ as the horocycle at $+\infty$).
We set $\hor(\om_1^k) = k \in \overline\Z$. The new boundary is
$\bd^* \T_q \cup \{ \om_1^k : k \in \overline\Z \}$. We write
$\wt\T_q$ for the new compactification, which we call the
\emph{horocyclic compactification.} It is induced by the metric
\begin{equation}\label{thetah}
\theta_h(z_1,w_1) = \theta(z_1,w_1) + 
\left| 
\frac{\hor(z_1)}{1+|\hor(z_1)|} - \frac{\hor(w_1)}{1+|\hor(w_1)|}
\right|\,,
\end{equation}
where $\theta$ is as in \eqref{theta} and we set 
$\frac{\pm\infty}{1+\infty}=\pm 1$. In this metric,
a sequence $(y_1^{(n)})$ tends to $\xi_1 \in \bd^*\T_q$ if and only if it
converges to $\xi_1$ in the end topology. It tends to $\om_1^k$ if
and only if $\hor(y_1^{(n)}) \to k$ ($k \in \overline\Z$) and $y_1^{(n)} \to \om_1$
in the end topology.

Again, we can take the closure $\wt\DL(q,r)$ of $\DL(q,r)$ in 
$\wt\T_q \times\wt\T_r$, the \emph{horocyclic compactification} of $\DL$.
In this case, the boundary consists of the following
$5$ disjoint pieces:
\begin{equation}\label{bdry-II}
\begin{gathered}
\bigl(\bd^*\T_q \times \{ \om_2^{-\infty} \}\bigr) \cup
\bigl(\{ \om_1^{-\infty} \} \times \bd^*\T_r\bigr) \cup
\bigl\{\om_1^{k}\om_2^{-k} : k \in \overline\Z\bigr\}\\
\cup
\bigl\{y_1\om_2^{-\hor(y_1)} : y_1 \in \T_q \bigr\} \cup
\bigl\{\om_1^{-\hor(y_2)}y_2 : y_2 \in \T_r \bigr\} \,.
\end{gathered}
\end{equation}
We omit the detailed description of convergence, which is a straightforward
adaptation of \eqref{convergence-I}. The mapping $\om_i^k \mapsto \om_i$
($i=1,2$) extends to a continuous surjection from the horocyclic onto the
geometric compactification, which restricted to $\DL(q,r)$ is the identity.

\begin{thm}\label{non-centred}
If $\al \ne 1/2$ then the Martin compactification of $\DL(q,r)$
with respect to $P = P_{\al}$ is the horocyclic compactification
$\wt{\DL}(q,r)$. The extension of the Martin kernel on the boundary
described in \eqref{bdry-II} is given by
\begin{alignat}{2}
K(x_1x_2, \xi_1\om_2^{-\infty}) &= K_1(x_1,\xi_1)\,,&\quad &\xi_1 \in
\bd^*\T_q\,,\tag{\rm i}\\[2pt]
K(x_1x_2, \om_1^{-\infty}\xi_2) &= K_2(x_2,\xi_2)\,,
&\quad &\xi_2\in \bd^*\T_r\,,\tag{\rm ii}\\[2pt]
K(x_1x_2, \om_1^{k}\om_2^{-k})  &=
\frac{\be^k + \be^{\hor(x_1)}}{\be^k + 1}\,,&\quad &k \in \overline\Z\,,
\tag{\rm iii}\\[2pt]
K(x_1x_2, y_1\om_2^{-\hor(y_1)}) &= \frac{\phi_{\al,q}(x_1,y_1)}
{\phi_{\al,q}(o_1,y_1)}\,,&\quad& y_1\in \T_q \,, \AND\tag{\rm iv}\\[2pt]
K(x_1x_2, \om_1^{-\hor(y_2)}y_2) &= \frac{\phi_{1-\al,r}(x_2,y_2)}
{\phi_{1-\al,r}(o_2,y_2)}\,,&\quad& y_2 \in \T_r\,.\tag{\rm v}
\end{alignat}
In \emph{(iii)}, $\be = (1-\al)/\al$, and for $k = \pm\infty$,
the right hand side is to be understood as the respective limit.
 
Each of the kernels in\/ {\rm (i)} and\/ {\rm (ii)} 
constitutes a minimal harmonic
function, while the ones of\/ {\rm (iii), (iv)} and\/ {\rm (v)} are non-minimal
harmonic. 
\end{thm}

\begin{proof} Once more, the proof that the minimal harmonic functions are precisley those
in  (i) and (ii) can be found in \cite{Woe}.

We now study convergence of $K(x,y)$ as $y$ tends to a boundary point.
This time, we use part (a) of Theorem \ref{estimates}.
We assume that $\al < 1/2$, since the case $\al > 1/2$ follows
by exchanging $q \leftrightarrow r$ and using $\al^* = 1-\al$ in the place of
$\al$, or also by using the relation \eqref{exchange}.
\smallskip

(i) Suppose that $y=y_1y_2 \to \xi_1\om_2^{-\infty}$ in $\wt{\DL}$.
Then $y_1 \to \xi_1$ and $y_2 \to \om_2$ in the end compactifications
of the respective trees. We proceed as in the proof of Theorem \ref{centred}
and find that in the formula of Theorem \ref{estimates}(a), the dominant 
one among the four terms in the $(...)$ on the right hand side is the
second one. It behaves like $\be^{-\dn_2}$. Therefore,  using 
\eqref{G1xy},
$$
G(x,y) \sim \frac{A_{\be}}{(q\,\be)^{\dn_1}\,(r\,\be)^{\dn_2}}
= A_{\be} \, F_1(x_1,y_1) \, F_2^*(x_2,y_2)\,,
$$
where (recall) $F_2^*(x_2,y_2)$ corresponds to exchanging 
$\al \leftrightarrow 1-\al$, that is, to the projection $P_{\al,r}$
onto $\T_r$ of $P_{\al}^* = P_{1-\al}$.
Therefore we obtain
$$
K(o,y) \sim K_1(x_1,y_1)\,K_2^*(x_2,y_2) \to K_1(x_1,\xi_1)\,K_2^*(x_2,\om_2)\,.
$$
Noting that $K_2^*(x_2,\om_2) = 1$, we get the proposed Martin kernel.

\smallskip

(ii) Similarly, if $y=y_1y_2 \to \om_1^{\infty}\xi_2$ then the dominant term
in the formula of Theorem \ref{estimates}(a) is the third one, which behaves
like $\be^{-\up_2}$. This and \eqref{G1xy} yield
$$
G(x,y) \sim \frac{A_{\be}}{(q\,\be)^{\dn_1}\,r^{\dn_2}\,\be^{\up_2}}
= A_{\be} \, F_1(x_1,y_1) \, F_2(x_2,y_2)\,,
$$
this time without passing to $P^*$. The conclusion is now as in (i) above.

\smallskip

(iii) Let $y \to \om_1^k\om_2^{-k}$, so that $\up_i \to \infty$ ($i=1,2$).

\smallskip

(a) $k=+\infty\,$. Then $\dn_1 - \up_1 = \up_2 - \dn_2 \to \infty$, so that   
the dominant term in the formula of Theorem \ref{estimates}(a) is the 
second one, as in (i). We get the same estimate as in (i), but have to
replace $\xi_1$ with $\om_1$, i.e.,
$$
K(x,y) \to K_1(x_1,\om_1)\,K_2^*(x_2,\om_2) = 1\,.
$$

(b) $k=-\infty\,$. In this case, the dominant term and asymptotic
estimate of $G(x,y)$ are the same as in (ii), whence
$$
K(o,y) \to K_1(x_1,\om_1)\,K_2(x_2,\om_2) = \be^{-\hor(x_2)} 
= \be^{\hor(x_1)} \,,
$$
since $K_1(x_1,\om_1)=1$.

\smallskip

(c) $k \in \Z$, and $\hor(y_1)=k$. 
In this case, all of $\up_i, \dn_i$ ($i=1,2$) tend to 
$\infty\,$. Also $\dn_1 - \up_1 = \up_2 - \dn_2 =  k - \hor(x_1)$. Therefore, 
in the formula of Theorem \ref{estimates}(a),  
among the four terms in the $(...)$ the second and the third one
are of the same order and dominate the other two. We obtain
$$
G(x,y) \sim \frac{A_{\be}}{(q\,\be)^{\dn_1}\,r^{\dn_2}}
\left(\frac{1}{\be^{\dn_2}} + \frac{1}{\be^{\up_2}}\right)
= A_{\be} \, F_1(x_1,y_1) \, F_2^*(x_2,y_2) \,
\bigl(1 + \be^{\hor(x_1)-k} \bigr)\,.
$$
Therefore
$$
K(x,y) \to K_1(x_1,\om_1)\,K_2^*(x_2,\om_2)\,
\frac{1 + \be^{\hor(x_1)-k}}{1 + \be^{-k}} = 
\frac{\be^k + \be^{\hor(x_1)}}{\be^k + 1}\,.
$$

\smallskip

(iv) Recall that when $y_1$ remains fixed and $y_2 \to \om_2^{-\hor(y_1)}$, 
then
$\dn_1 = \dn(x_1,y_1)$ and $\up_1 = \up(x_1,y_1)$ are constant, while
$\up_2 = \spn - \up_1$ and $\dn_2 = \spn-\dn_1$. 
In  Theorem \ref{estimates}(a), the first three of  
the four terms in the $(...)$ are of the same order and dominate fourth.
Thus
$$
G(x,y) \sim \frac{A_{\be}}{(q\,\be)^{\dn_1}\,r^{\dn_2}}\left(
\frac{B_{\be}}{\be^{\spn}} + \frac{\be^{\dn_1}-1}{\be^{\spn}} 
+ \frac{\be^{\up_1}-1}{\be^{\spn}}\right) =
A_{\be}  \,(\be-1)\,  F_2^*(x_2,y_2) \,\phi_{\al,q}(x_1,y_1)\,.
$$
This yields the proposed limit of $K(x,y)$.

\smallskip

(v) When $y_2$ is fixed and $y_1 \to \om_1^{-\hor(y_2)}$, we get analogously
$$
G(x,y) \sim \frac{A_{\be}}{(q\,\be)^{\dn_1}\,r^{\dn_2}}
\left(
\frac{B_{\be}}{\be^{\spn}} + \frac{\be^{\dn_1}-1}{\be^{\spn}} 
+ \frac{\be^{\up_1}-1}{\be^{\spn}}\right) =
A_{\be}  \,(\be-1)\,  F_1^*(x_1,y_1) \,\phi_{1-\al,r}(x_2,y_2)\,.
$$
Again, this yields the proposed limit of $K(x,y)$.
\end{proof}

\section{Positive eigenfunctions}\label{t-harmonic}
 
It is well known and easy to prove that 
positive \emph{$t$-harmonic} functions $h$ (satisfying $Ph = t\cdot h$) 
exist if and only if $t \ge \rho(P)$, see e.g. \cite{Wbook}, Lemma 7.2.
The Green kernel (resolvent) associated with eigenvalue $t$ is
$$
G(x,y|t) = \sum_{n=0}^{\infty} p^{(n)}(x,y)/t^n\,, \quad x,y \in X\,.
$$
(Instead of the variable $t$, often $z =1/t$ is used in the literature.)
The Martin compactification associated with $P$ and the eigenvalue $t$
can be constructed in the same way as described in the Introduction,
using the Martin kernel 
$$
K(x,y|t) = G(x,y|t)/G(o,y|t)\,.
$$
Now consider $P=P_{\al}$ on $\DL(q,r)$, its projections to the two trees,
and in particular, $\wt P$ on $\Z$. We fix 
$t \ge \rho= \rho(P) = 2\sqrt{\al(1-\al)}$.
Set
$$
\al(t) = \frac{t - \sqrt{t^2 - \rho^2}}{2t} \AND
\la(t) = \frac{t - \sqrt{t^2 - \rho^2}}{2\al}.
$$
Then the function on $\Z$ defined by $\psi(k) = \la(t)^k$ satisfies
$\wt P_{\al} \psi = t\cdot \psi$. We can lift this function to
$\T^1$, $\T_2$ and $\DL$ by using the respective projection, and we
obtain a $t$-harmonic function for the respective random walk.
Then we can conjugate the resepctive transtition matrix with the
lifted function, and divide by $t$. We end up with a new transition matrix.
On $\DL$, this becomes
\begin{equation}\label{conjugateP}
\frac{p_{\al}(x,y)\, \psi\bigl(\hor(y_1)\bigr)}{t \, \psi\bigl(\hor(x_1)\bigr)}
= p_{\al(t)}(x,y)\,.
\end{equation}
Consequently, the associated Green and Martin kernels on $\DL$ satisfy
\begin{equation}\label{conjugateG}
\begin{aligned}
G_{\al}(x,y|t) &= G_{\al(t)}(x,y)\, \la(t)^{\hor(x_1)-\hor(y_1)} \AND \\
K_{\al}(x,y|t) &= \frac{G_{\al}(x,y|t)}{G_{\al}(o,y|t)} = 
K_{\al(t)}(x,y)\, \la(t)^{\hor(x_1)}\,, 
\end{aligned}
\end{equation}
where $G_{\al(t)}(x,y)$ and $K_{\al(t)}(x,y)$ are the ordinary Green and Martin kernels 
(with $t=1$) of $P_{\al(t)}$ on $\DL$. Thus, the estimates of \S \ref{renewal}
also yield the asymptotics of $G_{\al}(x,y|t)$. Note here that 
$\al(\rho)=1/2$ and $\al(t) < 1/2$ when $t > \rho$.
Also note that formulas analogous to \eqref{conjugateP} and \eqref{conjugateG}
hold for the projected random walks on the two trees.

\begin{cor}\label{t-Martin}
The Martin compactification of $\DL(q,r)$ with respect to $P_{\al}$ and eigenvalue
$t$ is the geometric compactification $\wh{\DL}(q,r)$ when $t=\rho(P_{\al})$ and
the horocyclic compactification $\wt{\DL}(q,r)$ when $t>\rho(P_{\al})$. 
\end{cor}

We omit transcribing from \S \ref{martin} the explicit formulas for all the extended 
Martin kernels and just remark that for any $t \ge \rho$, we get
$$
\begin{aligned}
K(x_1x_2, \xi_1\om_2^{-\infty}|t) &= K_1(x_1,\xi_1|t)\,,
\quad \xi_1 \in \bd^*\T_q\,,\AND\\
K(x_1x_2, \om_1^{-\infty}\xi_2|t) &= K_2(x_1,\xi_1|t)\,,
\quad \xi_2 \in \bd^*\T_r\,.
\end{aligned}
$$
We have omitted the $\al$, resp. $1-\al$ in the subscripts, and 
the superscript of $\om_2^{\infty}$ ($i=1,2$) has to be omitted when $t=\rho$.

Once more, the Martin compactification is \emph{stable}
in the sense of Picardello and Woess \cite{PiWo}: in particular, 
the compactification
is the same for all $t > \rho$, while at the bottom of the positive spectrum, 
i.e., for $t=\rho$, it is smaller. Indeed, the identity on $\DL(q,r)$ extends 
to a continuous surjection from the horocyclic onto the geometric 
compactification.

\section{A remark on the elliptic Harnack inequality}\label{Harnack}

The \emph{elliptic Harnack inequality} for reversible random walks on
graphs appears frequently in recent research, see e.g. {\sc Hebisch
and Saloff-Coste} \cite{HeSC}, {\sc Delmotte} \cite{De}, 
{\sc Grigor'yan and Telcs} \cite{GrTe}, or -- most suitable in our context 
-- the recent note of {\sc Barlow} \cite{Bar}. Barlow shows among other that the elliptic Harnack
inequality for a random walk with ``controlled weights'' 
(in particular, for SRW)
on a graph $X$ with bounded vertex degrees) is equivalent with a Harnack
inequality for restricted Green functions
$$
G^D(x,y) = \sum_{n=0}^{\infty} \Prb[Z_n=y, Z_k \in D (k \le n) | Z_0=x]\,,
$$
where $D \subset X$ is finite.

In the formulation of \cite{Bar}, Theorem 2, this inequality -- denoted (HG) --
requires that there is a constant $C$ such that if $x_0, x, y \in X$
are such that $d(x_0,x)=d(x_0,y)=R \ge 1$ and $v \in D$ for all $v$ with
$d(x_0,v) \le 2R$, then
\begin{equation}\label{HG}
G^D(x_0,y) \le C\cdot G^D(x_0,x)\,.
\end{equation}
When the random walk is transient then we can let $D$ tend to $X$ (i.e.,
we use an increasing sequence $(D_n)$ of finite subsets whose union is $X$),
and we see that (HG) implies 
\begin{equation}\label{HGreen}
G(x_0,y) \le C \cdot G(x_0,x) \quad \text{for all}\; x_0, x, y \in X
\; \text{with}\; d(x_0,x)=d(x_0,y)\,.
\end{equation}
In \cite{Bar}, it is shown that the random walk on the lamplighter group which
corresponds to SRW on $\DL(2,2)$ does \emph{not} satisfy (HG), or equivalently,
the elliptic Harnack inequality.

This can also be seen easily from our asymptotic estimate. 
Indeed, consider SRW on $\DL(q,q)$ and $R \ge 1$.
We choose $x=x_1x_2$ such that $\hor(x_1)= d(o,x) =2R$, so that the 
relative position of $x$ with respect to $o$ is that of (III) in Figure 4. 
Also, we choose $y=y_1y_2$ such that $\hor(y_1)=0$ and $d(o,y) = 2R$,
with relative position as in (I) of Figure 4. Then, using Corollaries
\ref{II-III} and \ref{I} with $q=r$ and $\be=1$, we get
$$
G(o,x) \sim \frac{G_1(o_1,o_1)\,G_2(o_2,o_2)}{4\,R\, q^{2R}}
\AND
G(o,y) \sim \dfrac{q+1}{q-1}\,\frac{G_1(o_1,o_1)\,G_2(o_2,o_2)}{2\,R^2\,q^{R}}
\,,
$$
as $R \to \infty\,$. Thus, $G(o,x)/G(o,y) \to 0$, and 
\eqref{HGreen} does not hold.

\end{document}